\documentclass[sn-mathphys-num]{sn-jnl}

\usepackage{amsmath,amsfonts,amssymb}
\usepackage{enumerate}
\usepackage{graphics}
\usepackage{graphicx}
\usepackage{subcaption}
\usepackage{epsfig}
\usepackage{hhline}
\usepackage[boxruled,linesnumbered,procnumbered]{algorithm2e}
\usepackage{scalerel,stackengine}
\usepackage{ntheorem}
\usepackage{accents}
\usepackage[acronym]{glossaries}
\usepackage{booktabs,tabularx}
\stackMath
\newcommand\reallywidecheck[1]{%
\savestack{\tmpbox}{\stretchto{%
  \scaleto{%
    \scalerel*[\widthof{\ensuremath{#1}}]{\kern-.6pt\bigwedge\kern-.6pt}%
    {\rule[-\textheight/2]{1ex}{\textheight}}
  }{\textheight}%
}{0.5ex}}%
\stackon[1pt]{#1}{\scalebox{-1}{\tmpbox}}%
}

\usepackage[colorinlistoftodos]{todonotes}

\renewcommand{\Re}{{\mathbb R}}
\newcommand{\xbold}{{\bf x}}
\newcommand{\zbold}{{\bf z}}

\newcommand{\ubold}{{\bf u}}

\newcommand{\ybold}{{\bf y}}
\newcommand{\wbold}{{\bf w}}

\newcommand{\mc}[1]{\mathcal{#1}}

\newcommand{\R}{{\mathbb R}}
\renewcommand{\Re}{{\mathbb R}}
\newcommand{\Z}{{\mathbb Z}}

\newcounter{algo}

\newtheorem{theorem}{Theorem}
\newtheorem{example}{Example}%
\newtheorem{definition}{Definition}%
\newtheorem{proposition}{Proposition}%
\newtheorem{assumption}{Assumption}

\makeglossaries
\newacronym{wrt}{w.r.t.}{with respect to}
\newacronym{MINLP}{MINLP}{mixed-integer nonlinear problem}
\newacronym{MI-NEP}{MI-NEP}{mixed-integer Nash equilibrium problem}
\newacronym{NEP}{NEP}{Nash equilibrium problem}
\newacronym{BR}{BR}{best-response}
\newacronym{MI-GNEP}{MI-NEP}{mixed-integer generalized Nash equilibrium problem}
\newacronym{SOC}{SOC}{state of charge}
\newacronym{DSO}{DSO}{distribution system operator}
\newacronym{MVAD}{MVAD}{multi-vehicle automated driving}
\newacronym{MILP}{MILP}{Mixed-Integer Linear Program}

\newglossaryentry{MI-NE}
{
	name={MI-NE},
	description={mixed-integer Nash equilibrium},
	first={\glsentrydesc{MI-NE} (\glsentrytext{MI-NE})},
	plural={MI-NE},
	descriptionplural={mixed-integer Nash equilibria},
	firstplural={\glsentrydescplural{MI-NE} (\glsentryplural{MI-NE})}
}

\graphicspath{{figures/}}

\begin{document}

\title{Best-response algorithms for a class of monotone Nash equilibrium problems with mixed-integer variables}

\author[1]{\fnm{Filippo} \sur{Fabiani}}\email{filippo.fabiani@imtlucca.it}

\author*[2]{\fnm{Simone} \sur{Sagratella}}\email{sagratella@diag.uniroma1.it}

\affil*[1]{IMT School for Advanced Studies Lucca, Piazza San Francesco 19, 55100 Lucca, Italy}

\affil[2]{Department of Computer, Control and Management Engineering Antonio Ruberti, Sapienza University of Rome, Via Ariosto 25, 00185 Roma, Italy}

\abstract{We characterize the convergence properties of traditional \gls{BR} algorithms in computing solutions to \glspl{MI-NEP} that turn into a class of monotone \glspl{NEP} once relaxed the integer restrictions. We show that the sequence produced by a Jacobi/Gauss-Seidel \gls{BR} method always approaches a bounded region containing the entire solution set of the \gls{MI-NEP}, whose tightness depends on the problem data, and it is related to the degree of strong monotonicity of the relaxed \gls{NEP}. When the underlying algorithm is applied to the relaxed \gls{NEP}, we establish data-dependent complexity results characterizing its convergence to the unique solution of the \gls{NEP}. In addition, we derive one of the very few sufficient conditions for the existence of solutions to \glspl{MI-NEP}. The theoretical results developed bring important practical benefits, illustrated on a numerical instance of a smart building control application.}




\maketitle

\section{Introduction}\label{sec: introduction}
\glsresetall
In many real-world applications different decision makers frequently interact in a non-cooperative fashion to take optimal decisions that may depend on the opponents' strategies. If these decisions have to be made simultaneously, every agent is rational and has complete information about the other agents' optimization problems, then Nash equilibria can be considered as solutions of the resulting non-cooperative game. Theory and algorithms to solve \glspl{NEP} have been widely investigated in the literature, but mostly considering continuous strategy spaces for the agents -- see, e.g., \cite{FacchKanSu,scutari2012monotone,scutari2014real,grammatico2017dynamic,dreves2011solution,dreves2020nonsingularity,bigi2023approximate}. 
Despite in many crucial situations it is rather natural to restrict some variables to exclusively assume integer values, 
relevant contributions along this direction appeared only very recently in the literature and, to the best of our knowledge, they can be summarized as follows:
\begin{itemize}
\item Enumerative techniques for general instances. A branch-and-prune method to compute the whole set of solutions in the integer case was proposed in the seminal work \cite{sagratella2016computing}, which has successively been extended in \cite{sagratella2019generalized} to solve \glspl{MI-NEP} with linear coupling constraints. Very recently, the method has been further boosted to encompass non-linearities and non-convexities \cite{schwarze2023branch}.
\item Best-response algorithms for potential problems. In \cite{sagratella2017algorithms} has been proven that Gauss-Seidel \gls{BR} algorithms always converge to (approximate) \glspl{MI-NE} and there are no solutions that can not be computed in this way.
More recently, \cite{fabiani2022proximal} discussed how integer-compatible regularization functions enable for convergence either to an exact or approximate \gls{MI-NE}.

\item Best-response algorithms for 2-groups partitionable instances. This class of problems was defined in \cite{sagratella2016computing} and arises in several economics and engineering applications, see e.g. \cite{sagratella2017computing,passacantando2023finite}. Jacobi-type \gls{BR} algorithms can be effectively used to compute equilibria for these games.

 \item Convexification techniques for quasi-linear problems. These are problems in which, for fixed strategies of the opponent players, the cost function of every agent is linear and the respective strategy space is polyhedral. The method is based on the Nikaido–Isoda function and is able to recast the Nash game as a standard non-linear optimization problem, see \cite{harks2024generalized}.
\end{itemize}
Finally, we mention tailored solution algorithms for (potential) \glspl{MI-NEP} proposed for specific engineering applications spanning from automated driving \cite{fabiani2018mixed,fabiani2019multi} and smart mobility \cite{cenedese2021highway} to demand-side management \cite{cenedese2019charging}.

So far, there are no works on \glspl{MI-NEP} in which monotonicity of the relaxed, i.e., fully continuous, version associated has been explicitly used to study problem solvability and convergence of \gls{BR} methods. This is a fundamental gap in the literature, as monotonicity is a key condition to prove many theoretical properties for continuous \glspl{NEP} \cite{FacchPangBk}. 
We thus concentrate on the class of continuous \glspl{NEP} with contraction properties (see Assumption \ref{as: block-contraction}) introduced in \cite{scutari2012monotone}, which we extend to a mixed-integer setting by requiring additional conditions on the distance between continuous and mixed-integer \glspl{BR} (see Assumption~\ref{as: discrete bound}), and show that the resulting \glspl{MI-NEP} enjoy several interesting properties. 
Our main contributions can hence be listed as follows.
\begin{itemize}
 \item Focusing on a Jacobi/Gauss-Seidel \gls{BR} method (Algorithm \ref{alg: Jacobi}), we show that the sequence produced always approaches a bounded region containing the entire solution set of the \gls{MI-NEP} -- see Theorems \ref{th: error bound}, \ref{th: error bound wrt continuous NEP}. 
 The smaller the constants in Assumptions \ref{as: block-contraction} and \ref{as: discrete bound}, the smaller the size of this region. This brings several numerical advantages, as discussed in Section~\ref{sec: numerical}.
 Furthermore, in Section~\ref{sec: inexact} we also generalize Algorithm \ref{alg: Jacobi} by considering the computation of approximate mixed-integer \glspl{BR}, showing that similar theoretical results can be obtained.
\item We establish a relation between the contraction property defined in Assumption~\ref{as: block-contraction} and the degree of (strong) monotonicity of the problem -- see Section \ref{sec: assumption 1}. Specifically, we propose different ways to obtain Assumption \ref{as: block-contraction} by suitably perturbing any strongly monotone problem, see Proposition~\ref{pr: strong monotonicity}.

In addition, we consider several Nash problem structures to obtain different constant values satisfying Assumption \ref{as: discrete bound} -- Section \ref{sec: assumption 2}. 
As emphasized in those sections, the class of problems considered in this paper is wide and can thus be employed to model many real-world applications.
\item We consider the Jacobi-type continuous \gls{BR} method (Algorithm \ref{alg: Jacobi continuous}) to confirm convergence to the unique solution of the continuous \gls{NEP} relaxation of the original \gls{MI-NEP} (as stated in \cite{scutari2012monotone}).
In addition, i) we provide some complexity results -- see Theorem \ref{th: error bound continuous in continuous problems}, and ii) we show that Algorithm \ref{alg: Jacobi continuous} has similar convergence properties as Algorithm \ref{alg: Jacobi} for the \gls{MI-NEP} -- see Theorem \ref{th: error bound continuous}, also discussing the related numerical advantages in Section \ref{sec: numerical}.
\item We establish one of the very few sufficient conditions for the existence of solutions to \glspl{MI-NEP} -- see Proposition \ref{pr: existence}. Our conditions depend only on the constants defined in Assumptions \ref{as: block-contraction} and \ref{as: discrete bound}, and on the unique solution of the continuous \gls{NEP} relaxation of the original \gls{MI-NEP}. Therefore, they are general and can be used to define solvable problems. Moreover, if these conditions hold true, then there exists a unique \gls{MI-NE} and it can be computed through Algorithm~\ref{alg: Jacobi}.
\item In Section \ref{sec: numerical} we finally corroborate our theoretical findings  on a novel application involving the smart control of a building in which a number of residential units, endowed with some storing capacity, is interested in designing an optimal schedule to switch on/off high power domestic appliances over a prescribed time window to make the energy supply of the entire building smart and efficient. 
\end{itemize}

\section{Problem definition and main assumptions}\label{sec: preliminaries}
Consider a \gls{NEP} with $N$ players, indexed by the set $\mc I \triangleq \{1,\ldots,N\}$, and let $\nu \in \mc I$ be a generic player of the game.
We denote by
$x^\nu \in \mathbb R^{n_\nu}$ the vector representing the private strategies of the $\nu$-th player, and by $\xbold^{-\nu} \triangleq (x^{\nu'})_{\nu'\in\mc I \setminus\{\nu\}}$ the vector of all the other players strategies. We write $\mathbb R^{n}\ni\xbold \triangleq
(x^\nu,\xbold^{-\nu})$, where $n \triangleq \sum_{\nu \in \mc I} n_{\nu}$, to indicate the collective vector of strategies.
Any player $\nu$ has to solve an optimization problem that is parametric \gls{wrt} the other players variables, and the resulting \gls{NEP} reads as:
\begin{equation} \label{eq: prob}
\forall \nu \in \mc I:\left\{
    \begin{aligned}
        &\underset{x^\nu}{\min} && \theta_\nu (x^\nu, \xbold^{-\nu})\\
        &\textrm{ s.t. } && x^\nu \in \Omega_\nu \triangleq \{x^\nu \in X_\nu \mid x^\nu_j \in \Z, \, j = 1, \ldots, i_\nu\},
    \end{aligned}
\right.
\end{equation}
where the cost functions $\theta_\nu \, :\, \Re^{n} \, \rightarrow \, \Re$ are continuously differentiable and convex \gls{wrt} $x^\nu$, $X_\nu \subseteq \Re^{n_\nu}$ are (possibly unbounded) convex and closed sets, $i_\nu \leq n_\nu$ are nonnegative integers, and the feasible regions $\Omega_\nu$ are nonempty.
We will assume throughout the manuscript that any optimization problem $\nu$ in \eqref{eq: prob} always admits a solution for every given $\xbold^{-\nu}$.
We remark that if $i_\nu = 0$ for all $\nu$, i.e. the component-wise integer restrictions imposed through $x^\nu_j \in \Z, \, j = 1, \ldots, i_\nu$ vanish, then the collection of optimization problems in \eqref{eq: prob} is actually a classical \gls{NEP}, while in case $i_\nu > 0$, the Nash problem is a \gls{MI-NEP}.

We indicate the overall feasible set with $\Omega \triangleq \prod_{\nu \in \mc I} \Omega_\nu$, and its continuous relaxation with $X \triangleq \prod_{\nu \in \mc I} X_\nu$.
Let us introduce the \gls{BR} set for player $\nu$ at $\overline \xbold^{-\nu} \in \Omega_{-\nu} \triangleq \prod_{\nu' \in \mc I \setminus \{\nu\}} \Omega_{\nu'}$ as follows:
\begin{equation}\label{eq: best response}
 R_\nu(\overline \xbold^{-\nu}) \triangleq \arg \min_{x^\nu}~\theta_\nu (x^\nu, \overline \xbold^{-\nu}) \ \text{ s.t. } x^\nu \in \Omega_\nu.
\end{equation}
Note that $R_\nu(\overline \xbold^{-\nu})$ is nonempty for any $\overline \xbold^{-\nu} \in \Omega_{-\nu}$.
Computing an element of the \gls{BR} set requires, in general, the solution of a \gls{MINLP} -- see, e.g., \cite{belotti2013mixed,belotti2009branching,nowak2006relaxation,tawarmalani2002convexification}. In this work, we are interested in the following, standard notion of equilibrium for a Nash game.

\begin{definition}
A collective vector of strategies $\overline \xbold \in \Omega$ is a \gls{MI-NE} of the \gls{MI-NEP} in \eqref{eq: prob} if, for all $\nu \in \mc I$, it holds that
\begin{equation}\label{eq: optimality prob}
\theta_\nu(\overline x^\nu, \overline \xbold^{-\nu}) - \theta_\nu(\widehat x^\nu, \overline \xbold^{-\nu}) \leq 0, \, \text{ with } \widehat x^\nu \in R_\nu(\overline \xbold^{-\nu}).
\end{equation}
\hfill$\square$
\end{definition}

\noindent
Let us define the continuous \gls{BR} set for player $\nu$ at $\overline \xbold^{-\nu} \in \Omega_{-\nu}$:
\begin{equation}\label{eq: best response continuous}
 T_\nu(\overline \xbold^{-\nu}) \triangleq \arg \min_{x^\nu}~\theta_\nu (x^\nu, \overline \xbold^{-\nu}) \ \text{ s.t. } x^\nu \in X_\nu.
\end{equation}
For simplicity we assume that $T(\xbold) \triangleq \left( T_\nu(\xbold^{-\nu}) \right)_{\nu\in\mc I}$ is a point-to-point operator, which happens to be true if each $x^\nu\mapsto\theta_\nu (x^\nu, \xbold^{-\nu})$ in \eqref{eq: prob} is strictly convex.

The following assumptions define the class of problems we will deal with. Additional details about these assumptions, such as sufficient conditions to guarantee them, can be found in Section \ref{sec: assumptions}.
The first assumption is key to obtain interesting properties for the \gls{BR} algorithms we propose. As described in Section \ref{sec: assumptions}, such an assumption is related to certain monotonicity properties for the \gls{MI-NEP}, see Proposition \ref{pr: strong monotonicity}.
We denote with $\|\cdot\|$ any given norm.

\begin{assumption}\label{as: block-contraction}
There exist $\alpha \in [0,1)$ and $\wbold \in \Re^N$, $w_\nu>0$ for all $\nu \in \mc I$, such that, for all $\zbold, \ybold \in X$, the players' continuous \gls{BR} operators satisfy the following contraction property:
$$
\displaystyle \max_{\nu \in \mc I}~w_\nu \|T_\nu(\zbold^{-\nu}) - T_\nu(\ybold^{-\nu})\| \leq \alpha~\max_{\nu \in \mc I}~w_\nu \|z^\nu - y^\nu\|.
$$
\hfill$\square$
\end{assumption}
The function
$
\|\xbold\|^{\cal{B}(\wbold)} \triangleq \max_{\nu \in \mc I}~w_\nu \|x^\nu\|,
$
used in Assumption \ref{as: block-contraction}, is a norm if $\wbold \in \Re^N_{++}$ -- see Proposition \ref{pr: block norn} in Appendix.
If Assumption \ref{as: block-contraction} holds true, then $T$ is a block-contraction operator with modulus $\alpha$ and weight vector $\wbold$, and it is equivalent to writing $\|T(\zbold) - T(\ybold)\|^{\cal{B}(\wbold)} \leq \alpha \|\zbold - \ybold\|^{\cal{B}(\wbold)}$, see \cite[Sect. 3.1.2]{BertTsit97}.

We postulate next an assumption to bound the maximal distance between the set of mixed-integer \glspl{BR} and the continuous \gls{BR}, for all the players.

\begin{assumption}\label{as: discrete bound}
There exists $\beta > 0$ such that, for all $\nu \in \mc I$ and $\zbold^{-\nu} \in \Omega_{-\nu}$:
$$
\|\widehat x^\nu - T_\nu(\zbold^{-\nu})\| \leq \beta \sqrt{i_\nu}, \ \text{ for all } \, \widehat x^\nu \in R_\nu(\zbold^{-\nu}).
$$
\hfill$\square$
\end{assumption}

\noindent
In Section \ref{sec: assumptions} we show that classes of problems exist such that both the assumptions above hold true. A discussion about existence of solutions for this class of problems can be found, instead, in Section \ref{sec: existence}.

\section{Best-response methods}\label{sec: algorithms}

We focus on \gls{BR} methods for \glspl{MI-NEP}, as for instance Algorithm \ref{alg: Jacobi}, that is a general Jacobi-type method that incorporates different classical \gls{BR} algorithms. Depending on the sequence of indices sets $\{{\cal J}^k\}_{k \in \mathbb{N}}$ it can turn into, e.g., a Gauss-Seidel (sequential) algorithm, or a pure (parallel) Jacobi one.

\begin{algorithm}
\SetAlgoLined
Choose a starting point $\xbold^0 \in \Omega$ and set $k:=0$\;
\For{$k = 0, 1, \ldots$}{
 Select a subset ${\cal J}^k \subseteq \mc I$ of the players' indices\;
 \ForAll {$\nu \in {\cal J}^k$} {
  Compute a \gls{BR} ${\widehat x}^{k,\nu} \in R_\nu(\xbold^{k,-\nu})$\;
  Set $x^{k+1,\nu} := \widehat x^{k,\nu}$\;
 }
 \ForAll {$\nu \notin {\cal J}^k$} {
  Set $x^{k+1,\nu} := x^{k,\nu}$\;
 }
}
\caption{Jacobi-type method\label{alg: Jacobi}}
\end{algorithm}

\noindent
We study next the convergence properties of Algorithm \ref{alg: Jacobi} under Assumptions \ref{as: block-contraction} and \ref{as: discrete bound}:

\begin{theorem}\label{th: error bound}
Suppose that Assumptions~\ref{as: block-contraction} and \ref{as: discrete bound} hold true, and that
\begin{equation}\label{eq: i > 0}
\max_{\nu \in \mc I}~\sqrt{i_\nu} > 0.
\end{equation}
In addition, assume that, in Algorithm \ref{alg: Jacobi}, every $h$ iterations at least one \gls{BR} of any player $\nu$ is computed,
that is, $\nu \in \cup_{t = k}^{k + h} {\cal J}^t$ for each player $\nu$ and each iterate $k$.
Let $\{\xbold^k\}_{k \in \mathbb{N}} \subseteq \Omega$ be the sequence generated by Algorithm \ref{alg: Jacobi}, and let ${\cal S}$ be the (possibly empty) set of the equilibria of the \gls{MI-NEP} defined in \eqref{eq: prob}.
\begin{enumerate}[(i)]
\item For every $\gamma \in \left(1, \tfrac{1}{\alpha} \right)$ and every $\xbold^* \in {\cal S}$, Algorithm~\ref{alg: Jacobi} generates a point $\xbold^{k}$ such that
\begin{equation}\label{eq: error bound}
 \displaystyle \max_{\nu \in \mc I}~w_\nu \|x^{k,\nu} - x^{*,\nu}\| < 2 \beta \left( \tfrac{\gamma}{1 - \gamma \alpha} \right) \max_{\nu \in \mc I}~w_\nu \sqrt{i_\nu},
\end{equation}
after at most $\overline k_\gamma$ iterations, with
$$
\overline k_\gamma \triangleq h \Bigl\lceil \log_\gamma \left( \max \left\{ \tfrac{(1- \gamma \alpha)~\max_{\zbold^* \in {\cal S}, \nu \in \mc I}~w_\nu \|x^{0,\nu} - z^{*,\nu}\|}{2 \beta \max_{\nu \in \mc I}~w_\nu \sqrt{i_\nu}}, \gamma \right\} \right) \Bigr\rceil.
$$
\item For for every $\gamma \in \left(1, \tfrac{1}{\alpha} \right)$ and every $\xbold^* \in {\cal S}$, any point $\xbold^{k}$, with $k \geq \overline k_\gamma$, satisfies \eqref{eq: error bound}.
\item Assume $\Omega$ bounded or ${\cal S}$ nonempty. Every cluster point $\widetilde \xbold$ of $\{\xbold^k\}_{k \in \mathbb{N}}$ (at least one exists) is contained in $\Omega$ and satisfies the following inequality for every $\xbold^* \in {\cal S}$:
\begin{equation}\label{eq: error bound lim}
 \displaystyle \max_{\nu \in \mc I}~w_\nu \|\widetilde x^{\nu} - x^{*,\nu}\| \leq 2 \beta \left( \tfrac{1}{1 - \alpha} \right)~\max_{\nu \in \mc I}~w_\nu \sqrt{i_\nu}.
\end{equation} 
\end{enumerate}
\hfill$\square$
\end{theorem}
\textit{Proof.}
(i) Assume without loss of generality that all $\xbold^k$, with $k < \overline k_\gamma$, violate \eqref{eq: error bound} for some $\xbold^* \in {\cal S}$.
For any $k < \overline k_\gamma$, the following chain of inequalities holds:
\begin{align*}
 \displaystyle 
 & \left( \tfrac{1}{\gamma} \right) \max_{\nu \in \mc I}~w_\nu \|x^{k,\nu} - x^{*,\nu}\| \overset{\text{(a)}}{\geq} \alpha \max_{\nu \in \mc I}~w_\nu \|x^{k,\nu} - x^{*,\nu}\| + 2 \beta \max_{\nu \in \mc I}~w_\nu \sqrt{i_\nu} \\
 & \qquad \overset{\text{(b)}}{\geq} \max_{\nu \in \mc I}~w_\nu \|T_\nu(\xbold^{k, -\nu}) - T_\nu(\xbold^{*, -\nu})\| + 2 \beta \max_{\nu \in \mc I}~w_\nu \sqrt{i_\nu} \\
 & \qquad \overset{\text{(c)}}{\geq} \max_{\nu \in \mc I}~w_\nu \|T_\nu(\xbold^{k, -\nu}) - T_\nu(\xbold^{*, -\nu})\| + \max_{\nu \in {\cal J}^{k}}~w_\nu \|T_\nu(\xbold^{k, -\nu}) - x^{k+1, \nu}\|\\ 
 & \hspace{7cm}  + \max_{\nu \in {\cal J}^{k}}~w_\nu \|T_\nu(\xbold^{*, -\nu}) - x^{*, \nu}\| \\
 & \qquad \overset{\text{(d)}}{\geq} \max_{\nu \in {\cal J}^{k}}~w_\nu \|x^{k+1, \nu} - x^{*, \nu}\|,
\end{align*}
where (a) is a direct consequence of the fact that $\xbold^k$ violates \eqref{eq: error bound} for $\xbold^*$, while (b) and (c) follows by Assumption~\ref{as: block-contraction} and \ref{as: discrete bound}, respectively. Finally, (d) is a consequence of the following observation: let $\overline \nu \in {\cal J}^{k}$ be a player such that $w_{\overline \nu}~\|x^{k+1, \overline \nu} - x^{*, \overline \nu}\| = \max_{\nu \in {\cal J}^{k}}~w_\nu \|x^{k+1, \nu} - x^{*, \nu}\|$, then we get
\begin{align*}
 \displaystyle 
 &  \max_{\nu \in \mc I}~w_\nu \|T_\nu(\xbold^{k, -\nu}) - T_\nu(\xbold^{*, -\nu})\| + \\ 
 &  \;\;\;\;\; \max_{\nu \in {\cal J}^{k}}~w_\nu \|T_\nu(\xbold^{k, -\nu}) - x^{k+1, \nu}\| + \max_{\nu \in {\cal J}^{k}}~w_\nu \|T_\nu(\xbold^{*, -\nu}) - x^{*, \nu}\| \\
 & \geq w_{\overline \nu}~\|T_{\overline \nu}(\xbold^{k, -{\overline \nu}}) - T_{\overline \nu}(\xbold^{*, -{\overline \nu}})\| + w_{\overline \nu}~\|x^{k+1, {\overline \nu}} - T_{\overline \nu}(\xbold^{k, -{\overline \nu}})\|+\\
 &\;\;\;\;\; w_{\overline \nu}~\|T_{\overline \nu}(\xbold^{*, -{\overline \nu}}) - x^{*, {\overline \nu}}\| \geq w_{\overline \nu}~\|x^{k+1, \overline \nu} - x^{*, \overline \nu}\|.
\end{align*}
Therefore, reasoning by induction, we can conclude that
$
\max_{\nu \in \mc I}~w_\nu \|x^{0,\nu} - x^{*,\nu}\| \geq \cdots \geq \max_{\nu \in \mc I}~w_\nu \|x^{\overline k_\gamma,\nu} - x^{*,\nu}\|,
$
and, by the definition of $h$, for any $k \leq \overline k_\gamma - h$, we obtain
$
(1/\gamma) \max_{\nu \in \mc I}~w_\nu \|x^{k,\nu} - x^{*,\nu}\| \geq \max_{\nu \in \mc I}~w_\nu \|x^{k+h,\nu} - x^{*,\nu}\|.
$
As a consequence, we obtain directly that
\begin{equation}\label{eq: decrease norm j}
\left( \tfrac{1}{\gamma} \right)^{\left( \tfrac{\overline k_\gamma}{h} \right)} \max_{\nu \in \mc I}~w_\nu \|x^{0,\nu} - x^{*,\nu}\| \geq \max_{\nu \in \mc I}~w_\nu \|x^{\overline k_\gamma, \nu} - x^{*,\nu}\|.
\end{equation}
Now, observe that, since $x^0$ violates \eqref{eq: error bound} for $\xbold^*$, we have
$$
\tfrac{(1- \gamma \alpha) \max_{\zbold^* \in {\cal S}, \nu \in \mc I}~w_\nu \|x^{0,\nu} - z^{*,\nu}\|}{2 \beta \max_{\nu \in \mc I}~w_\nu \sqrt{i_\nu}} \geq \gamma,
$$
and hence
$
\overline k_\gamma = h \Bigl\lceil \log_\gamma \left( \tfrac{(1- \gamma \alpha) \max_{\zbold^* \in {\cal S}, \nu \in \mc I}~w_\nu \|x^{0,\nu} - z^{*,\nu}\|}{2 \beta \max_{\nu \in \mc I}~w_\nu \sqrt{i_\nu}} \right) \Bigr\rceil.
$
This implies
$$
\gamma^{\left( \tfrac{\overline k_\gamma}{h} \right)} \geq \tfrac{(1 - \gamma \alpha) \max_{\nu \in \mc I}~w_\nu \|x^{0,\nu} - x^{*,\nu}\|}{2 \beta \max_{\nu \in \mc I}~w_\nu \sqrt{i_\nu}},
$$
and since $\gamma > 1$,
$
 2 \beta \left( \tfrac{\gamma}{1 - \gamma \alpha} \right) \max_{\nu \in \mc I}~w_\nu \sqrt{i_\nu} > \left( \tfrac{1}{\gamma} \right)^{\left( \tfrac{\overline k_\gamma}{h} \right)} \max_{\nu \in \mc I}~w_\nu \|x^{0,\nu} - x^{*,\nu}\|.
$
This latter relation combined with \eqref{eq: decrease norm j} shows that $\xbold^{\overline k_\gamma}$ meets \eqref{eq: error bound} for $\xbold^*$.

(ii) By relying on (i), without loss of generality we can assume that $\xbold^{\overline k_\gamma}$ satisfies \eqref{eq: error bound} for some $\xbold^* \in {\cal S}$.
If $\arg \max_{\nu \in \mc I}~w_\nu \|x^{\overline k_\gamma + 1,\nu} - x^{*,\nu}\| \cap {\cal J}^{\overline k_\gamma} = \emptyset$,
then we trivially obtain $\max_{\nu \in \mc I}~w_\nu \|x^{\overline k_\gamma + 1,\nu} - x^{*,\nu}\| \leq \max_{\nu \in \mc I}~w_\nu \|x^{\overline k_\gamma,\nu} - x^{*,\nu}\|$ and, therefore, $\xbold^{\overline k_\gamma + 1}$ satisfies \eqref{eq: error bound} for $\xbold^*$. Alternatively, it holds that:
\begin{align*}
 \displaystyle 
 & \max_{\nu \in \mc I}~w_\nu \|x^{\overline k_\gamma + 1,\nu} - x^{*,\nu}\| = \max_{\nu \in {\cal J}^{\overline k_\gamma}}~w_\nu \|x^{\overline k_\gamma + 1,\nu} - x^{*,\nu}\| \\
 & \qquad \leq \max_{\nu \in {\cal J}^{\overline k_\gamma}}~w_\nu \|T_\nu(\xbold^{\overline k_\gamma, -\nu}) - T_\nu(\xbold^{*, -\nu})\| + \\ 
 & \qquad \;\;\;\;\; \max_{\nu \in {\cal J}^{\overline k_\gamma}}~w_\nu \|T_\nu(\xbold^{\overline k_\gamma, -\nu}) - x^{\overline k_\gamma+1, \nu}\| + \max_{\nu \in {\cal J}^{\overline k_\gamma}}~w_\nu \|T_\nu(\xbold^{*, -\nu}) - x^{*, \nu}\| \\
 & \qquad \overset{\text{(a)}}{\leq} \max_{\nu \in {\cal J}^{\overline k_\gamma}}~w_\nu \|T_\nu(\xbold^{\overline k_\gamma, -\nu}) - T_\nu(\xbold^{*, -\nu})\| + 2 \beta \max_{\nu \in \mc I}~w_\nu \sqrt{i_\nu} \\
 & \qquad \overset{\text{(b)}}{\leq} \alpha \max_{\nu \in \mc I}~w_\nu \|x^{\overline k_\gamma,\nu} - x^{*,\nu}\| + 2 \beta \max_{\nu \in \mc I}~w_\nu \sqrt{i_\nu} \\
 & \qquad \overset{\text{(c)}}{<} 2 \beta \left( \tfrac{\gamma \alpha}{1 - \gamma \alpha} + 1 \right) \max_{\nu \in \mc I}~w_\nu \sqrt{i_\nu} = 2 \beta \left( \tfrac{1}{1 - \gamma \alpha} \right) \max_{\nu \in \mc I}~w_\nu \sqrt{i_\nu},
\end{align*}
where (a) and (b) follow by Assumption~\ref{as: discrete bound} and \ref{as: block-contraction}, respectively, while (c) since $\xbold^{\overline k_\gamma}$ satisfies \eqref{eq: error bound} for $\xbold^*$. Then, also in this case $\xbold^{\overline k_\gamma + 1}$ satisfies \eqref{eq: error bound} for $\xbold^*$, as $\gamma > 1$. By iterating this reasoning, we can conclude that any point $\xbold^{k}$, with $k \geq \overline k_\gamma$, satisfies \eqref{eq: error bound} for $\xbold^*$, and hence for any \gls{MI-NE} of the \gls{MI-NEP} in \eqref{eq: prob}.

(iii) Observe that, in view of (ii), the sequence $\{\xbold^k\}_{k \in \mathbb{N}}$ is entirely contained in a bounded set for all $k \geq \overline k_\gamma$. Therefore, at least one cluster point $\widetilde \xbold$ shall exist. Since $\Omega$ is closed, we obtain $\widetilde \xbold \in \Omega$. Following the same steps as in the proof of part (ii) and recalling the definition of $h$, the following inequality holds true for every $k \geq \overline k_\gamma$ and $\xbold^* \in {\cal S}$, and for some $1 \leq j \leq h$:
$$
 \displaystyle 
 \max_{\nu \in \mc I}~w_\nu \|x^{k+j,\nu} - x^{*,\nu}\| \leq \alpha \max_{\nu \in \mc I}~w_\nu \|x^{k,\nu} - x^{*,\nu}\| + 2 \beta \max_{\nu \in \mc I}~w_\nu \sqrt{i_\nu}.
$$
Therefore, we obtain that
\begin{align*}
 \displaystyle \max_{\nu \in \mc I}~w_\nu \|\widetilde x^{\nu} - x^{*,\nu}\| & \leq 
 \left( \lim_{t \to +\infty} \alpha^t \right) 2 \beta \left( \tfrac{\gamma}{1 - \gamma \alpha} \right) \max_{\nu \in \mc I}~w_\nu \sqrt{i_\nu} + \\
 & \quad \left( \sum_{t = 0}^{+\infty} \alpha^t \right) 2 \beta \max_{\nu \in \mc I}~w_\nu \sqrt{i_\nu} = 2 \beta \left( \tfrac{1}{1 - \alpha} \right) \max_{\nu \in \mc I}~w_\nu \sqrt{i_\nu},
\end{align*}
where the equality above is valid because $\alpha < 1$. That is, the thesis is true.
\hfill$\blacksquare$

\noindent
Theorem \ref{th: error bound} characterizes the sequence produced by Algorithm \ref{alg: Jacobi}. In \eqref{eq: error bound} a bound for the maximal distance expressed in terms of the norm $\|\cdot\|^{\cal{B}(\wbold)}$ between the sequence generated by the algorithm and every solution of the \gls{MI-NEP} is defined. This bound defines a region that contains all the solutions of the \gls{MI-NEP} defined in \eqref{eq: prob} and it is strictly related to the values of $\alpha$ and $\beta$ defined in Assumptions \ref{as: block-contraction} and \ref{as: discrete bound}, respectively. Specifically the bound decreases, and as a consequence the corresponding convergence region shrinks, if $\alpha$ or $\beta$ descrease.
With item (i) of Theorem \ref{th: error bound} we ensure that the sequence produced by the algorithm reaches this convergence region in no more than $\overline k_\gamma$ iterations, and with item (ii) we know that the sequence remains in this convergence region, once reached.
Item (iii) of Theorem \ref{th: error bound} considers cluster points of the sequence and provides a slightly refined bound defined in \eqref{eq: error bound lim}. Note that this nice behavior for the sequence produced by Algorithm \ref{alg: Jacobi} of being attracted by the convergence region defined by the bound in \eqref{eq: error bound} is guaranteed also if the $\Omega$ is unbounded and at least one solution of the \gls{MI-NEP} exists.

\sloppy For all $\nu \in \mc I$, we define quantity
$
D_\nu \triangleq \max_{z^\nu, y^\nu \in X_\nu}~\|z^\nu - y^\nu\|.
$
If every set $X_\nu$ is bounded, then every $D_\nu$ is finite, and hence the bound $\overline k_\gamma$, defined in item (i) of Theorem \ref{th: error bound}, explicitly reads as:
$
\overline k_\gamma = h \Bigl\lceil \log_\gamma \left( \max \left\{ \tfrac{(1- \gamma \alpha) \max_{\nu \in \mc I}~w_\nu D_\nu}{2 \beta \max_{\nu \in \mc I}~w_\nu \sqrt{i_\nu}}, \gamma \right\} \right) \Bigr\rceil.
$
The following example shows that, under Assumptions~\ref{as: block-contraction} and \ref{as: discrete bound}, convergence results stronger than those in Theorem \ref{th: error bound} are not possible. Specifically, the sequence produced by Algorithm \ref{alg: Jacobi} may not converge to the unique solution of the \gls{MI-NEP}, but it can be attracted by the convergence region containing the solution, according to Theorem \ref{th: error bound}.

\begin{example}\label{ex: no convergence}
 Consider a \gls{MI-NEP} with 2 players, each controlling a single integer variable $i_1 = i_2 = 1$, cost functions and private constraints as:
 $\theta_1(x^1,x^2) = (x^1)^2 + (1+\varepsilon) x^1 x^2$, $X_1 = [l^1,u^1]$, with $l^1 \leq -1$ and $u^1 \geq 1$, 
 $\theta_2(x^1,x^2) = (x^2)^2 - (1+\varepsilon) x^1 x^2$, $X_2 = [l^2,u^2]$, with $l^2 \leq -1$ and $u^2 \geq 1$, and $\varepsilon > 0$ small.
 The unique equilibrium $\xbold^*$ of this problem is the origin, i.e., $\xbold^*=(0,0)^{\scriptscriptstyle\top}$.
 
 Consider now the sequence produced by Algorithm \ref{alg: Jacobi} starting from $\xbold^1 = (-1,1)^{\scriptscriptstyle\top}$. Player 1 is at an optimal point, while player 2 moves to $\xbold^2 = (-1,-1)^{\scriptscriptstyle\top}$. Now is player 1 the one that is not at an optimal point and moves to $\xbold^3 = (1,-1)^{\scriptscriptstyle\top}$. Again, player 2 moves and player 1 is at an optimal point: $\xbold^4 = (1,1)^{\scriptscriptstyle\top}$. Finally, the movement of player 1, with player 2 at an optimal point, brings the iterations back to the starting point: $\xbold^5 = (-1,1)^{\scriptscriptstyle\top} = \xbold^1$.
 Notice that this sequence is actually independent from the choice of the players included in ${\cal J}^k$ at each iteration $k$, and it is not convergent. None of the 4 cluster points $\{\xbold^1, \xbold^2, \xbold^3, \xbold^4\}$ is a solution of the Nash problem.

 Depending on $\varepsilon$, the \gls{MI-NEP} meets Assumptions \ref{as: block-contraction} and \ref{as: discrete bound}.
 Considering the Euclidean norm, we obtain $\alpha = (1+\varepsilon)/2$ and $\wbold = (1,1)^{\scriptscriptstyle\top}$ (see Proposition \ref{pr: sufficient for assumption 1} in Section \ref{sec: assumptions}) and $\beta = 1/2$ (see Proposition \ref{pr: sufficient for beta 2} in Section \ref{sec: assumptions}).
 Therefore, as shown in item (iii) of Theorem \ref{th: error bound}, none of the cluster point is far from the solution more than a given bound:
 $
 \max_{\nu \in \{1, 2\}}~\|x^{j,\nu} - x^{*,\nu}\|_2 \leq 2(1 - \varepsilon)^{-1}, \, j = 1, 2, 3, 4.
 $
Referring to Theorem~\ref{th: error bound}.(i) and (ii), a similar distance from the solution is finally obtained in less than $\overline k_\gamma$ iterations from any other starting point.
\hfill$\square$
\end{example}

\subsection{About continuous NEPs}

Let us consider the fully continuous case, i.e., \eqref{eq: i > 0} is not verified and then $\max_{\nu \in \mc I}~\sqrt{i_\nu} = 0$.
In this case Algorithm \ref{alg: Jacobi} is actually equivalent to the continuous version of the Jacobi-type \gls{BR} method defined in Algorithm~\ref{alg: Jacobi continuous}. 

\begin{algorithm}[h!]
\SetAlgoLined
Choose a starting point $\xbold^0 \in X$ and set $k:=0$\;
\For{$k = 0, 1, \ldots$}{
 Select a subset ${\cal J}^k \subseteq \mc I$ of the players' indices\;
 \ForAll {$\nu \in {\cal J}^k$} {
  Compute the continuous \gls{BR} ${\widehat x}^{k,\nu} = T_\nu(\xbold^{k,-\nu})$\;
  Set $x^{k+1,\nu} := \widehat x^{k,\nu}$\;
 }
 \ForAll {$\nu \notin {\cal J}^k$} {
  Set $x^{k+1,\nu} := x^{k,\nu}$\;
 }
}
\caption{Jacobi-type continuous method\label{alg: Jacobi continuous}}
\end{algorithm}

The following result shows that, under Assumption~\ref{as: block-contraction} (notice that Assumption \ref{as: discrete bound} is meaningless in the continuous setting), Algorithm \ref{alg: Jacobi continuous} converges to the unique equilibrium of the \gls{NEP} with continuous variables.

\begin{theorem}\label{th: error bound continuous in continuous problems}
Suppose that Assumption~\ref{as: block-contraction} holds true, and that $\max_{\nu \in \mc I}~\sqrt{i_\nu} = 0$.
Assume that, in Algorithm \ref{alg: Jacobi continuous}, every $h$ iterations at least one \gls{BR} of any player $\nu$ is computed,
that is, $\nu \in \cup_{t = k}^{k + h} {\cal J}^t$ for each player $\nu$ and each iterate $k$.
Let $\{\xbold^k\}_{k \in \mathbb{N}} \subseteq \Omega$ be the sequence generated by Algorithm \ref{alg: Jacobi continuous}. Problem \eqref{eq: prob} has a unique equilibrium $\overline \xbold$, and the following statements hold true:
\begin{enumerate}[(i)]
\item Assume $\xbold^0 \neq \overline \xbold$ and, for every $\varepsilon > 0$, let
$$
\widehat k_\varepsilon \triangleq h \Bigl\lceil \log_{\alpha} \left( \min \left\{ \tfrac{\varepsilon}{\max_{\nu \in \mc I}~w_\nu \|x^{0,\nu} - \overline x^{\nu}\|}, 1 \right\} \right) \Bigr\rceil.
$$
For all $k \geq \widehat k_\varepsilon$, every point $\xbold^{k}$ satisfies
\begin{equation}\label{eq: error bound continuous in continuous problems}
 \displaystyle \max_{\nu \in \mc I}~w_\nu \|x^{k,\nu} - \overline x^{\nu}\| \leq \varepsilon.
\end{equation}
\item The sequence $\{\xbold^k\}_{k \in \mathbb{N}}$ converges to the unique equilibrium $\overline \xbold$ of the \gls{NEP}.
\end{enumerate}
\hfill$\square$
\end{theorem}
\textit{Proof.}
Existence of a solution is guaranteed by standard results, see, e.g., \cite{FacchPangBk}. About uniqueness, assume by contradiction that the \gls{NEP} has a solution $\widetilde \xbold \neq \overline \xbold$. By using Assumption \ref{as: block-contraction}, we obtain $
\alpha \max_{\nu \in \mc I}~w_\nu \|\overline x^\nu - \widetilde x^\nu\| \geq 
\max_{\nu \in \mc I}~w_\nu \|T_\nu(\overline \xbold^{-\nu}) - T_\nu(\widetilde \xbold^{-\nu})\| 
= \max_{\nu \in \mc I}~w_\nu \|\overline x^\nu - \widetilde x^\nu\|,
$
which is incompatible with $\alpha < 1$. Thus, the \gls{NEP} has a unique solution.

For all $k > 0$, we have: $\max_{\nu \in \mc I}~w_\nu \|x^{k,\nu} - \overline x^{\nu}\| = \max\left\{ \max_{\nu \in {\cal J}^{k-1}}~w_\nu \|T_\nu(\xbold^{k-1, -\nu}) - T_\nu(\overline \xbold^{-\nu})\|, \max_{\nu \not\in {\cal J}^{k-1}}~w_\nu \|x^{k-1,\nu} - \overline x^{\nu}\| \right\}$ $\leq \max\left\{ \alpha \max_{\nu \in {\cal J}^{k-1}}~w_\nu \|x^{k-1,\nu} - \overline x^{\nu}\|, \max_{\nu \not\in {\cal J}^{k-1}}~w_\nu \|x^{k-1,\nu} - \overline x^{\nu}\| \right\},$ 
where the inequality is due to Assumption~\ref{as: block-contraction}.
By the definition of $h$, for every $k \geq h$ we then obtain:
\begin{equation}\label{eq: continuous chain}
 \max_{\nu \in \mc I}~w_\nu \|x^{k,\nu} - \overline x^{\nu}\| \leq \alpha \max_{\nu \in \mc I}~w_\nu \|x^{k-h,\nu} - \overline x^{\nu}\|.
\end{equation}
(i) Assume by contradiction that $\xbold^k$, with $k \geq \widehat k_\varepsilon$, violates \eqref{eq: error bound continuous in continuous problems}. In this case we obtain the following chain of inequalities that can not be verified:
\begin{equation*}
 \displaystyle 
 \varepsilon < \max_{\nu \in \mc I}~w_\nu \|x^{k,\nu} - \overline x^{\nu}\|
 \overset{\text{(a)}}{\leq} \alpha^{\tfrac{\widehat k_\varepsilon}{h}} \max_{\nu \in \mc I}~w_\nu \|x^{0,\nu} - \overline x^{\nu}\| \overset{\text{(b)}}{\leq} \varepsilon,
\end{equation*}
where (a) comes from \eqref{eq: continuous chain} and (b) is due to the definition of $k_\varepsilon$.

(ii) In view of \eqref{eq: continuous chain}, it holds that
$
 \lim_{k \to \infty} \max_{\nu \in \mc I}~w_\nu \|x^{k,\nu} - \overline x^{\nu}\|
 \leq \lim_{k \to \infty} \alpha^{k/h} \max_{\nu \in \mc I}~w_\nu \|x^{0,\nu} - \overline x^{\nu}\| = 0.
$
The thesis hence follows since the function $\|\cdot\|^{\cal{B}(\wbold)}$ is a norm -- see Proposition \ref{pr: block norn} in Appendix.
\hfill$\blacksquare$

\noindent
A similar convergence result to Theorem \ref{th: error bound continuous in continuous problems} (ii) is stated in \cite{scutari2012monotone} and it can be established from \cite[Prop.~1.1 in \S 3.1.1, Prop.~1.4 in \S 3.1.2]{BertTsit97} if one considers only the Gauss-Seidel version of Algorithm \ref{alg: Jacobi continuous}. 
Note that the complexity measure in item (i) of Theorem \ref{th: error bound continuous in continuous problems} is original and it is useful to predict the number of iterations needed to compute an approximate solution through Algorithm \ref{alg: Jacobi continuous}.

\subsection{Relations between MI-NEPs and their continuous NEP relaxations, and a discussion about existence of solutions}\label{sec: existence}

We consider \glspl{MI-NEP} satisfying Assumptions \ref{as: block-contraction} and \ref{as: discrete bound} and such that \eqref{eq: i > 0} is verified to make an analysis on the sequence produced by Algorithm \ref{alg: Jacobi continuous} \gls{wrt} the solution set of the \gls{MI-NEP}. To complete the picture, we also consider the sequence produced by Algorithm \ref{alg: Jacobi} and define relations with the unique solution of the continuous \gls{NEP} relaxation of the original \gls{MI-NEP}.  

The following result provides better bounds for Algorithm \ref{alg: Jacobi continuous} than those established in Theorem \ref{th: error bound} for Algorithm \ref{alg: Jacobi}.


\begin{theorem}\label{th: error bound continuous}
Suppose that Assumptions~\ref{as: block-contraction} and \ref{as: discrete bound} hold true, and that \eqref{eq: i > 0} is verified.
Assume that, in Algorithm \ref{alg: Jacobi continuous}, every $h$ iterations at least one \gls{BR} of any player $\nu$ is computed,
that is, $\nu \in \cup_{t = k}^{k + h} {\cal J}^t$ for each player $\nu$ and each iterate $k$.
Let $\{\xbold^k\}_{k \in \mathbb{N}} \subseteq X$ be the sequence generated by Algorithm \ref{alg: Jacobi continuous}, and let ${\cal S}$ be the (possibly empty) set of the equilibria of the \gls{MI-NEP} in \eqref{eq: prob}.
\begin{enumerate}[(i)]
\item For every $\gamma \in \left(1, \tfrac{1}{\alpha} \right)$ and every $\xbold^* \in {\cal S}$, Algorithm \ref{alg: Jacobi continuous} generates a point $\xbold^{k}$ such that
\begin{equation}\label{eq: error bound continuous}
 \displaystyle \max_{\nu \in \mc I}~w_\nu \|x^{k,\nu} - x^{*,\nu}\| < \beta \left( \tfrac{\gamma}{1 - \gamma \alpha} \right) \max_{\nu \in \mc I}~w_\nu \sqrt{i_\nu},
\end{equation}
after at most $\widetilde k_\gamma$ iterations, with
$$
\widetilde k_\gamma \triangleq h \Bigl\lceil \log_\gamma \left( \max \left\{ \tfrac{(1- \gamma \alpha) \max_{\zbold^* \in {\cal S}, \nu \in \mc I} w_\nu \|x^{0,\nu} - z^{*,\nu}\|}{\beta \max_{\nu \in \mc I} w_\nu \sqrt{i_\nu}}, \gamma \right\} \right) \Bigr\rceil.
$$
\item For every $\gamma \in \left(1, \tfrac{1}{\alpha} \right)$ and every $\xbold^* \in {\cal S}$, any point $\xbold^{k}$, with $k \geq \widetilde k_\gamma$, satisfies \eqref{eq: error bound continuous}.
\item The sequence $\{\xbold^k\}_{k \in \mathbb{N}}$ converges to a unique point $\overline \xbold \in X$. The following inequality holds for every $\xbold^* \in {\cal S}$:
\begin{equation}\label{eq: error bound lim continuous}
 \displaystyle \max_{\nu \in \mc I}~w_\nu \|\overline x^{\nu} - x^{*,\nu}\| \leq \beta \left( \tfrac{1}{1 - \alpha} \right) \max_{\nu \in \mc I}~w_\nu \sqrt{i_\nu}.
\end{equation} 
\end{enumerate}
\hfill$\square$
\end{theorem}
\textit{Proof.}
The proof is similar to that of Theorem \ref{th: error bound}.\\
(i) Assume without loss of generality that all $\xbold^k$, with $k < \widetilde k_\gamma$, violate \eqref{eq: error bound continuous} for some $\xbold^* \in {\cal S}$.
The following chain of inequalities holds for any $k < \widetilde k_\gamma$:
\begin{align*}
 \displaystyle 
 & \left( \tfrac{1}{\gamma} \right) \max_{\nu \in \mc I}~w_\nu \|x^{k,\nu} - x^{*,\nu}\| \\
 & \qquad \overset{\text{(a)}}{\geq} \alpha \max_{\nu \in \mc I}~w_\nu \|x^{k,\nu} - x^{*,\nu}\| + \beta \max_{\nu \in \mc I}~w_\nu \sqrt{i_\nu} \\
 & \qquad \overset{\text{(b)}}{\geq} \max_{\nu \in \mc I}~w_\nu \|T_\nu(\xbold^{k, -\nu}) - T_\nu(\xbold^{*, -\nu})\| + \beta \max_{\nu \in \mc I}~w_\nu \sqrt{i_\nu} \\
 & \qquad \overset{\text{(c)}}{\geq} \max_{\nu \in \mc I}~w_\nu \|T_\nu(\xbold^{k, -\nu}) - T_\nu(\xbold^{*, -\nu})\| + \max_{\nu \in {\cal J}^{k}}~w_\nu \|T_\nu(\xbold^{*, -\nu}) - x^{*, \nu}\| \\
 & \qquad \overset{\text{(d)}}{\geq} \max_{\nu \in {\cal J}^{k}}~w_\nu \|x^{k+1, \nu} - x^{*, \nu}\|,
\end{align*}
where (a) is a direct consequence of the fact that $\xbold^k$ violates \eqref{eq: error bound} for $\xbold^*$, while (b) and (c) follow by Assumption~\ref{as: block-contraction} and \ref{as: discrete bound}, respectively. Inequality (d), instead, is a consequence of the following observation: let $\overline \nu \in {\cal J}^{k}$ be a player such that $w_{\overline \nu}~\|x^{k+1, \overline \nu} - x^{*, \overline \nu}\| = \max_{\nu \in {\cal J}^{k}}~w_\nu \|x^{k+1, \nu} - x^{*, \nu}\|$. Then, we obtain: $ \max_{\nu \in \mc I}~w_\nu \|T_\nu(\xbold^{k, -\nu}) - T_\nu(\xbold^{*, -\nu})\| + \max_{\nu \in {\cal J}^{k}}~w_\nu \|T_\nu(\xbold^{*, -\nu}) - x^{*, \nu}\| \geq w_{\overline \nu}~\|T_{\overline \nu}(\xbold^{k, -{\overline \nu}}) - T_{\overline \nu}(\xbold^{*, -{\overline \nu}})\| + w_{\overline \nu}~\|T_{\overline \nu}(\xbold^{*, -{\overline \nu}}) - x^{*, {\overline \nu}}\|  = w_{\overline \nu}~\|x^{k+1, \overline \nu} - T_{\overline \nu}(\xbold^{*, -{\overline \nu}})\| + w_{\overline \nu}~\|T_{\overline \nu}(\xbold^{*, -{\overline \nu}}) - x^{*, {\overline \nu}}\|\geq w_{\overline \nu}~\|x^{k+1, \overline \nu} - x^{*, \overline \nu}\|$.
By following the same reasoning as in the proof of Theorem \ref{th: error bound}.(i), we thus have
\begin{equation}\label{eq: decrease norm j continuous}
\left( \tfrac{1}{\gamma} \right)^{\left( \tfrac{\widetilde k_\gamma}{h} \right)} \max_{\nu \in \mc I}~w_\nu \|x^{0,\nu} - x^{*,\nu}\| \geq \max_{\nu \in \mc I}~w_\nu \|x^{\widetilde k_\gamma, \nu} - x^{*,\nu}\|.
\end{equation}
Now observe that, since $x^0$ violates \eqref{eq: error bound continuous} for $\xbold^*$, then:
$$
\tfrac{(1- \gamma \alpha) \max_{\zbold^* \in {\cal S}, \nu \in \mc I}~w_\nu \|x^{0,\nu} - z^{*,\nu}\|}{\beta \max_{\nu \in \mc I}~w_\nu \sqrt{i_\nu}} \geq \gamma,
$$
and, as a consequence,
$$
\widetilde k_\gamma = h \Bigl\lceil \log_\gamma \left( \tfrac{(1- \gamma \alpha) \max_{\zbold^* \in {\cal S}, \nu \in \mc I}~w_\nu \|x^{0,\nu} - z^{*,\nu}\|}{\beta \max_{\nu \in \mc I}~w_\nu \sqrt{i_\nu}} \right) \Bigr\rceil,
$$
which in turn implies that
$$
\gamma^{\left( \tfrac{\widetilde k_\gamma}{h} \right)} \geq \tfrac{(1 - \gamma \alpha) \max_{\nu \in \mc I}~w_\nu \|x^{0,\nu} - x^{*,\nu}\|}{\beta \max_{\nu \in \mc I}~w_\nu \sqrt{i_\nu}},
$$
and therefore that
$$
 \beta \left( \tfrac{\gamma}{1 - \gamma \alpha} \right) \max_{\nu \in \mc I}~w_\nu \sqrt{i_\nu} > \left( \tfrac{1}{\gamma} \right)^{\left( \tfrac{\widetilde k_\gamma}{h} \right)} \max_{\nu \in \mc I}~w_\nu \|x^{0,\nu} - x^{*,\nu}\|.
$$
Combining this latter relation with \eqref{eq: decrease norm j continuous} shows that $\xbold^{\widetilde k_\gamma}$ satisfies \eqref{eq: error bound continuous} for $\xbold^*$.

(ii) By relying on item (i), without loss of generality we can assume that $\xbold^{\widetilde k_\gamma}$ satisfies \eqref{eq: error bound continuous} for some $\xbold^* \in {\cal S}$. If $\arg \max_{\nu \in \mc I}~w_\nu \|x^{\widetilde k_\gamma + 1,\nu} - x^{*,\nu}\| \cap {\cal J}^{\widetilde k_\gamma} = \emptyset$,
then we trivially obtain that $\max_{\nu \in \mc I}~w_\nu \|x^{\widetilde k_\gamma + 1,\nu} - x^{*,\nu}\| \leq \max_{\nu \in \mc I}~w_\nu \|x^{\widetilde k_\gamma,\nu} - x^{*,\nu}\|$, and therefore $\xbold^{\widetilde k_\gamma + 1}$ satisfies \eqref{eq: error bound continuous} for $\xbold^*$.
Alternatively, the following chain of inequalities holds true:
\begin{align*}
 \displaystyle 
 & \max_{\nu \in \mc I}~w_\nu \|x^{\widetilde k_\gamma + 1,\nu} - x^{*,\nu}\| = \max_{\nu \in {\cal J}^{\widetilde k_\gamma}}~w_\nu \|x^{\widetilde k_\gamma + 1,\nu} - x^{*,\nu}\|\\
 & \qquad = \max_{\nu \in {\cal J}^{\widetilde k_\gamma}}~w_\nu \|T_\nu(\xbold^{\widetilde k_\gamma, -\nu}) - x^{*,\nu}\| \\
 & \qquad \leq \max_{\nu \in {\cal J}^{\widetilde k_\gamma}}~w_\nu \|T_\nu(\xbold^{\widetilde k_\gamma, -\nu}) - T_\nu(\xbold^{*, -\nu})\| + \max_{\nu \in {\cal J}^{\widetilde k_\gamma}}~w_\nu \|T_\nu(\xbold^{*, -\nu}) - x^{*, \nu}\| \\
 & \qquad \overset{\text{(a)}}{\leq} \max_{\nu \in {\cal J}^{\widetilde k_\gamma}}~w_\nu \|T_\nu(\xbold^{\widetilde k_\gamma, -\nu}) - T_\nu(\xbold^{*, -\nu})\| + \beta \max_{\nu \in \mc I}~w_\nu \sqrt{i_\nu} \\
 & \qquad \overset{\text{(b)}}{\leq} \alpha \max_{\nu \in \mc I}~w_\nu \|x^{\widetilde k_\gamma,\nu} - x^{*,\nu}\| + \beta \max_{\nu \in \mc I}~w_\nu \sqrt{i_\nu} \\
 & \qquad \overset{\text{(c)}}{<} \beta \left( \tfrac{\gamma \alpha}{1 - \gamma \alpha} + 1 \right) \max_{\nu \in \mc I}~w_\nu \sqrt{i_\nu} = \beta \left( \tfrac{1}{1 - \gamma \alpha} \right) \max_{\nu \in \mc I}~w_\nu \sqrt{i_\nu},
\end{align*}
where (a) and (b) follow from Assumption~\ref{as: discrete bound} and \ref{as: block-contraction}, respectively, while (c) holds true since $\xbold^{\widetilde k_\gamma}$ satisfies \eqref{eq: error bound continuous} for $\xbold^*$. The proof hence follows as in Theorem \ref{th: error bound}.(ii).

(iii) Theorem~\ref{th: error bound continuous in continuous problems} guarantees the convergence to the unique point $\overline \xbold \in X$. Similar to the proof of item (ii), and recalling the definition of $h$, the following inequality holds true for every $k \geq \widetilde k_\gamma$ and $\xbold^* \in {\cal S}$, and for some $1 \leq j \leq h$:
$
 \max_{\nu \in \mc I}~w_\nu \|x^{k+j,\nu} - x^{*,\nu}\| \leq \alpha \max_{\nu \in \mc I}~w_\nu \|x^{k,\nu} - x^{*,\nu}\| + \beta \max_{\nu \in \mc I}~w_\nu \sqrt{i_\nu}.
$
We thus have that: $\max_{\nu \in \mc I}~w_\nu \|\overline x^{\nu} - x^{*,\nu}\|\leq 
\left( \lim_{t \to +\infty} \alpha^t \right) \beta \left( \tfrac{\gamma}{1 - \gamma \alpha} \right) \max_{\nu \in \mc I}~w_\nu \sqrt{i_\nu} +\left( \sum_{t = 0}^{+\infty} \alpha^t \right) \beta \max_{\nu \in \mc I}~w_\nu \sqrt{i_\nu}= \beta \left( \tfrac{1}{1 - \alpha} \right) \max_{\nu \in \mc I}~w_\nu \sqrt{i_\nu}$,
which proves the desired claim.
\hfill$\blacksquare$


\noindent
Comparing Theorem \ref{th: error bound continuous} and Theorem \ref{th: error bound}, we note that the bounds provided by \eqref{eq: error bound continuous} and \eqref{eq: error bound lim continuous} for the sequence produced by Algorithm \ref{alg: Jacobi continuous} are tighter than those shown in \eqref{eq: error bound} and \eqref{eq: error bound lim} for the sequence produced by Algorithm \ref{alg: Jacobi}.

Note that, by Theorem \ref{th: error bound continuous in continuous problems} we know that $\overline \xbold$ in Theorem~\ref{th: error bound continuous}.(iii) is actually the unique solution of the continuous \gls{NEP} relaxation of the \gls{MI-NEP}. Therefore, \eqref{eq: error bound lim continuous} provides also an upper bound for the distance between every \gls{MI-NE} of the \gls{MI-NEP} in \eqref{eq: prob} and the unique solution of its continuous \gls{NEP} relaxation.

The following result characterizes instead the distance of the sequence produced by Algorithm \ref{alg: Jacobi} from the unique solution of the continuous \gls{NEP} relaxation of the \gls{MI-NEP}. We omit the proof for the sake of presentation, as it is similar to those of Theorems \ref{th: error bound} and \ref{th: error bound continuous}.

\begin{theorem}\label{th: error bound wrt continuous NEP}
Suppose that Assumptions~\ref{as: block-contraction} and \ref{as: discrete bound} hold true, and that \eqref{eq: i > 0} is met.
Moreover, assume that, in Algorithm \ref{alg: Jacobi}, every $h$ iterations at least one \gls{BR} of any player $\nu$ is computed, i.e., $\nu \in \cup_{t = k}^{k + h} {\cal J}^t$ for each player $\nu$ and iterate $k$.
Let $\{\xbold^k\}_{k \in \mathbb{N}} \subseteq \Omega$ be the sequence generated by Algorithm \ref{alg: Jacobi}, and let $\overline \xbold$ be the unique solution of the continuous \gls{NEP} relaxation of the \gls{MI-NEP} in \eqref{eq: prob}.
\begin{enumerate}[(i)]
\item For every $\gamma \in \left(1, \tfrac{1}{\alpha} \right)$, Algorithm~\ref{alg: Jacobi} generates a point $\xbold^{k}$ such that
\begin{equation}\label{eq: error bound wrt continuous NEP}
 \displaystyle \max_{\nu \in \mc I}~w_\nu \|x^{k,\nu} - \overline x^{\nu}\| < \beta \left( \tfrac{\gamma}{1 - \gamma \alpha} \right) \max_{\nu \in \mc I}~w_\nu \sqrt{i_\nu},
\end{equation}
after at most $\widehat {\widehat k}_\gamma$ iterations, with
$$
\widehat {\widehat k}_\gamma \triangleq h \Bigl\lceil \log_\gamma \left( \max \left\{ \tfrac{(1- \gamma \alpha) \max_{\nu \in \mc I}~w_\nu \|x^{0,\nu} - \overline x^{\nu}\|}{\beta \max_{\nu \in \mc I}~w_\nu \sqrt{i_\nu}}, \gamma \right\} \right) \Bigr\rceil.
$$
\item For every $\gamma \in \left(1, \tfrac{1}{\alpha} \right)$, any point $\xbold^{k}$, with $k \geq \widehat {\widehat k}_\gamma$, satisfies \eqref{eq: error bound wrt continuous NEP}.
\item Every cluster point $\widetilde \xbold$ of $\{\xbold^k\}_{k \in \mathbb{N}}$ (at least one exists) is contained in $\Omega$ and satisfies the following inequality:
\begin{equation}\label{eq: error bound lim continuous wrt continuous NEP}
 \displaystyle \max_{\nu \in \mc I}~w_\nu \|\widetilde x^{\nu} - \overline x^{\nu}\| \leq \beta \left( \tfrac{1}{1 - \alpha} \right) \max_{\nu \in \mc I}~w_\nu \sqrt{i_\nu}.
\end{equation} 
\end{enumerate}
\hfill$\square$
\end{theorem}

\noindent
The results in Theorem \ref{th: error bound wrt continuous NEP} are valid even when the \gls{MI-NEP} in \eqref{eq: prob} does not admit any equilibrium or the feasible region $\Omega$ is unbounded, since in those cases the sequence generated by Algorithm~\ref{alg: Jacobi} would stay in a bounded set.

Finally, we note that Theorems \ref{th: error bound continuous} and \ref{th: error bound wrt continuous NEP} allow us to establish sufficient conditions for the existence of solutions to the \gls{MI-NEP}.

\begin{proposition}\label{pr: existence}
 Suppose that Assumptions~\ref{as: block-contraction} and \ref{as: discrete bound} hold true, and that \eqref{eq: i > 0} is verified.
 Let $\overline \xbold$ be the unique solution of the continuous \gls{NEP} relaxation of the \gls{MI-NEP} in \eqref{eq: prob}.
 If for every $\widetilde\xbold \in \Omega$ such that \eqref{eq: error bound lim continuous wrt continuous NEP} holds the integer components are unique, i.e., $\widetilde x^\nu_j = \overline x^\nu_j$, $j = 1, \ldots, i_\nu$, $\forall \nu \in {\mathcal I}$, for some $\overline \xbold \in \Omega$, then there exists a unique solution $\xbold^*$ of the \gls{MI-NEP} such that $x^{*,\nu}_j = \overline x^\nu_j$, $j = 1, \ldots, i_\nu$, $\forall \nu \in {\mathcal I}$, and the sequence produced by Algorithm \ref{alg: Jacobi} converges to $\xbold^*$.
\hfill$\square$
\end{proposition}
\textit{Proof.}
By Theorem \ref{th: error bound continuous}.(iii), every solution $\xbold^*$ of the \gls{MI-NEP} satisfies \eqref{eq: error bound lim continuous}. Therefore, $x^{*,\nu}_j = \overline x^\nu_j$, $j = 1, \ldots, i_\nu$, $\forall \nu \in {\mathcal I}$, and the set of solutions of the \gls{MI-NEP} contains at most one point, as otherwise Assumption \ref{as: block-contraction} is violated.

Theorem \ref{th: error bound wrt continuous NEP} and the assumptions made guarantee that, in finite iterations, each $\xbold^k$ produced by Algorithm \ref{alg: Jacobi} satisfies $x^{k,\nu}_j = \overline x^\nu_j$, $j = 1, \ldots, i_\nu$, $\forall \nu \in {\mathcal I}$. Thus, eventually, Algorithm \ref{alg: Jacobi} is equivalent to Algorithm \ref{alg: Jacobi continuous} on the continuous variables while the integer ones are fixed to the integer components of $\overline \xbold$. From Theorem \ref{th: error bound continuous in continuous problems}, Algorithm \ref{alg: Jacobi} converges to a solution of the \gls{MI-NEP}.
\hfill$\blacksquare$

\noindent
The following example describes a \gls{MI-NEP} satisfying the sufficient conditions given in Proposition \ref{pr: existence}.

\begin{example}\label{ex: existence}
 Consider a \gls{MI-NEP} with 2 players, each controlling a single integer variable $i_1 = i_2 = 1$, cost functions and private constraints as:
 $\theta_1(x^1,x^2) = \tfrac{1}{2}\left(x^1- \eta_1\right)^2 + \varepsilon_1 x^1 x^2$, $X_1 = [l^1,u^1]$, 
 $\theta_2(x^1,x^2) = \tfrac{1}{2}\left(x^2- \eta_2\right)^2 + \varepsilon_2 x^1 x^2$, $X_2 = [l^2,u^2]$, where $\eta \in X$ and $\varepsilon_1, \varepsilon_2$ are small numbers.

 The unique solution of the continuous \gls{NEP} relaxation is 
 $
 \overline \xbold=\left(\tfrac{\eta_1-\varepsilon_1 \eta_2}{1-\varepsilon_1 \varepsilon_2},\tfrac{\eta_2-\varepsilon_2 \eta_1}{1-\varepsilon_1 \varepsilon_2}\right)^{\scriptscriptstyle\top}.
 $
 Depending on $\varepsilon_1$ and $\varepsilon_2$, this \gls{MI-NEP} satisfies Assumptions \ref{as: block-contraction} and \ref{as: discrete bound}.
 Specifically, considering the Euclidean norm, we obtain $\alpha = \max\left\{|\varepsilon_1|,|\varepsilon_2|\right\}$ and $\wbold = (1,1)^{\scriptscriptstyle\top}$ (see Proposition \ref{pr: sufficient for assumption 1} in Section \ref{sec: assumptions}) and $\beta = 1/2$ (see Proposition \ref{pr: sufficient for beta 2} in Section \ref{sec: assumptions}).
 In this case condition \eqref{eq: error bound lim continuous wrt continuous NEP} reads as
 $$
  \displaystyle
  \left|\overline x^1 \!-\! \widetilde x^{1}\right|\leq \tfrac{1}{2} \left( \tfrac{1}{1 - \max\left\{|\varepsilon_1|,|\varepsilon_2|\right\}} \right), \; 
    \left|\overline x^2 \!-\! \widetilde x^{2}\right|\leq \tfrac{1}{2} \left( \tfrac{1}{1 - \max\left\{|\varepsilon_1|,|\varepsilon_2|\right\}} \right).
$$
 Let us assume without loss of generality that $\overline x^1 - \lfloor \overline x^1 \rfloor \leq \lceil \overline x^1 \rceil - \overline x^1$ and $\overline x^2 - \lfloor \overline x^2 \rfloor \leq \lceil \overline x^2 \rceil - \overline x^2$. If it holds that
 $$
   \lceil \overline x^1 \rceil> \overline x^1 + \tfrac{1}{2} \left( \tfrac{1}{1 - \max\left\{|\varepsilon_1|,|\varepsilon_2|\right\}} \right), \;
   \lceil \overline x^2 \rceil > \overline x^2 + \tfrac{1}{2} \left( \tfrac{1}{1 - \max\left\{|\varepsilon_1|,|\varepsilon_2|\right\}} \right),
$$
 then $\widetilde \xbold = \lfloor \overline \xbold \rfloor$ is the unique point in $\Omega$ that satisfies condition \eqref{eq: error bound lim continuous wrt continuous NEP}, and it is the unique solution to the \gls{MI-NEP} by Proposition \ref{pr: existence}.
 \hfill$\square$
\end{example}

\noindent
The following example illustrates that sufficient conditions for the existence of solutions based on the constants $\alpha$ and $\beta$ only may not be established.

\begin{example}\label{ex: no existence}
 Consider a \gls{MI-NEP} with 2 players, each controlling a single integer variable $i_1 = i_2 = 1$, cost functions and private constraints as:
 $\theta_1(x^1,x^2) = \tfrac{1}{2}\left(x^1 - \tfrac{1}{2}\right)^2 + \varepsilon x^1 \left(x^2 - \tfrac{1}{2}\right)$, $X_1 = [0,1]$, 
 $\theta_2(x^1,x^2) = \tfrac{1}{2}\left(x^2- \tfrac{1}{2}\right)^2 - \varepsilon \left(x^1 - \tfrac{1}{2}\right) x^2$, $X_2 = [0,1]$, where $\varepsilon$ is a small positive number.

  Depending on $\varepsilon$, this \gls{MI-NEP} satisfies Assumptions \ref{as: block-contraction} and \ref{as: discrete bound}.
 Specifically, with the Euclidean norm we obtain $\alpha = \varepsilon$ and $\wbold = (1,1)^{\scriptscriptstyle\top}$ (see Proposition \ref{pr: sufficient for assumption 1} in Section \ref{sec: assumptions}) and $\beta = 1/2$ (see Proposition \ref{pr: sufficient for beta 2} in Section \ref{sec: assumptions}). Therefore the value of $\alpha$ can be arbitrarily small, but this problem does not admit any solution for every value of $\varepsilon$:
 $\theta_1(0,0) = \tfrac{1}{8}$, $\theta_2(0,0) = \tfrac{1}{8}$,
 $\theta_1(1,0) = \tfrac{1}{8} - \tfrac{\varepsilon}{2}$, $\theta_2(1,0) = \tfrac{1}{8}$,
 $\theta_1(0,1) = \tfrac{1}{8}$, $\theta_2(0,1) = \tfrac{1}{8} + \tfrac{\varepsilon}{2}$,
 $\theta_1(1,1) = \tfrac{1}{8} + \tfrac{\varepsilon}{2}$, $\theta_2(1,1) = \tfrac{1}{8} - \tfrac{\varepsilon}{2}$.
 \hfill$\square$
\end{example}

\subsection{Remark about inexact BR algorithms} \label{sec: inexact}
Algorithm \ref{alg: Jacobi} requires the computation of an element of the \gls{BR} set, thereby calling for the solution of a \gls{MINLP} that may be prohibitive in most practical cases. For this reason, we consider now the realistic scenario in which approximate \glspl{BR} are computed in place of exact ones.
Specifically, we assume the existence of $\varepsilon \geq 0$ such that $\widehat x^{k,\nu}$ in step 5 of Algorithm \ref{alg: Jacobi} is computed to satisfying:
\begin{equation}\label{eq: inexact BR}
  \|\widehat x^{k,\nu} - \widetilde x^{k,\nu} \| \leq \varepsilon, 
\end{equation}
for some $\widetilde x^{k,\nu} \in R_\nu(\xbold^{k,-\nu})$.
In particular, we note that Theorem \ref{th: error bound} applies also in this case with the following modifications:
\begin{itemize}
    \item In item (i), \eqref{eq: error bound} is replaced with
    \begin{equation}\label{eq: error bound inexact}
 \displaystyle \max_{\nu \in \mc I}~w_\nu \|x^{k,\nu} - x^{*,\nu}\| < \left( \tfrac{\gamma}{1 - \gamma \alpha} \right) \max_{\nu \in \mc I}~w_\nu (2 \beta \sqrt{i_\nu} + \varepsilon),
\end{equation}
and
$$
\overline k_\gamma \triangleq h \Bigl\lceil \log_\gamma \left( \max \left\{ \tfrac{(1- \gamma \alpha)~\max_{\zbold^* \in {\cal S}, \nu \in \mc I}~w_\nu \|x^{0,\nu} - z^{*,\nu}\|}{\max_{\nu \in \mc I}~w_\nu (2 \beta \sqrt{i_\nu} + \varepsilon)}, \gamma \right\} \right) \Bigr\rceil;
$$

\item In item (ii), \eqref{eq: error bound} is replaced with \eqref{eq: error bound inexact};

\item In item (iii), \eqref{eq: error bound lim} is replaced with
\begin{equation}\label{eq: error bound lim inexact}
 \displaystyle \max_{\nu \in \mc I}~w_\nu \|\widetilde x^{\nu} - x^{*,\nu}\| \leq \left( \tfrac{1}{1 - \alpha} \right)~\max_{\nu \in \mc I}~w_\nu (2 \beta \sqrt{i_\nu} + \varepsilon);
\end{equation} 

\item In the proof of item (i):
\begin{itemize}
 \item Inequality
$$
\left( \tfrac{1}{\gamma} \right) \max_{\nu \in \mc I}~w_\nu \|x^{k,\nu} - x^{*,\nu}\| \geq \max_{\nu \in {\cal J}^{k}}~w_\nu \|x^{k+1, \nu} - x^{*, \nu}\|,
$$
still holds because
\begin{align*}
 \displaystyle 
 & \left( \tfrac{1}{\gamma} \right) \max_{\nu \in \mc I}~w_\nu \|x^{k,\nu} - x^{*,\nu}\| \\
 & \qquad \geq \alpha \max_{\nu \in \mc I}~w_\nu \|x^{k,\nu} - x^{*,\nu}\| + \max_{\nu \in \mc I}~w_\nu (2 \beta \sqrt{i_\nu} + \varepsilon) \\
 & \qquad \geq \max_{\nu \in \mc I}~w_\nu \|T_\nu(\xbold^{k, -\nu}) - T_\nu(\xbold^{*, -\nu})\| + \max_{\nu \in \mc I}~w_\nu (2 \beta \sqrt{i_\nu} + \varepsilon) \\
 & \qquad \geq \max_{\nu \in \mc I}~w_\nu \|T_\nu(\xbold^{k, -\nu}) - T_\nu(\xbold^{*, -\nu})\| + \\
 & \qquad\quad \max_{\nu \in {\cal J}^{k}}~w_\nu \|T_\nu(\xbold^{k, -\nu}) - \widetilde x^{k, \nu}\| + \\ 
 & \qquad\quad \max_{\nu \in {\cal J}^{k}}~w_\nu \|T_\nu(\xbold^{*, -\nu}) - x^{*, \nu}\| + \max_{\nu \in \mc I}~w_\nu \varepsilon \\
 & \qquad \geq \max_{\nu \in {\cal J}^{k}}~w_\nu \|\widetilde x^{k, \nu} - x^{*, \nu}\| + \max_{\nu \in \mc I}~w_\nu \varepsilon \\
 & \qquad \overset{(a)}{\geq} \max_{\nu \in {\cal J}^{k}}~w_\nu \|\widetilde x^{k, \nu} - x^{*, \nu}\| + \max_{\nu \in \mc I}~w_\nu \|x^{k+1,\nu} - \widetilde x^{k, \nu}\| \\
 & \qquad \geq \max_{\nu \in {\cal J}^{k}}~w_\nu (\|\widetilde x^{k, \nu} - x^{*, \nu}\| + \|x^{k+1,\nu} - \widetilde x^{k, \nu}\|) \\
 & \qquad \geq \max_{\nu \in {\cal J}^{k}}~w_\nu \|x^{k+1, \nu} - x^{*, \nu}\|,
\end{align*}
where $\widetilde x^{k,\nu}$ is defined in \eqref{eq: inexact BR} and the first four inequalities are similar to those in the proof of Theorem \ref{th: error bound}, and $(a)$ follows from \eqref{eq: inexact BR};

\item Relation in \eqref{eq: decrease norm j} is the same;

\item Moreover, it holds that
$$
\tfrac{(1- \gamma \alpha) \max_{\zbold^* \in {\cal S}, \nu \in \mc I}~w_\nu \|x^{0,\nu} - z^{*,\nu}\|}{\max_{\nu \in \mc I}~w_\nu (2 \beta \sqrt{i_\nu} + \varepsilon)} \geq \gamma,
$$
and hence
$
\overline k_\gamma = h \Bigl\lceil \log_\gamma \left( \tfrac{(1- \gamma \alpha) \max_{\zbold^* \in {\cal S}, \nu \in \mc I}~w_\nu \|x^{0,\nu} - z^{*,\nu}\|}{\max_{\nu \in \mc I}~w_\nu (2 \beta \sqrt{i_\nu} + \varepsilon)} \right) \Bigr\rceil;
$
thus
$$
\gamma^{\left( \tfrac{\overline k_\gamma}{h} \right)} \geq \tfrac{(1 - \gamma \alpha) \max_{\nu \in \mc I}~w_\nu \|x^{0,\nu} - x^{*,\nu}\|}{\max_{\nu \in \mc I}~w_\nu (2 \beta \sqrt{i_\nu} + \varepsilon)},
$$
and, by \eqref{eq: decrease norm j}:
\begin{align*}
 \left( \tfrac{\gamma}{1 - \gamma \alpha} \right) \max_{\nu \in \mc I}~w_\nu (2 \beta \sqrt{i_\nu} + \varepsilon) &> \left( \tfrac{1}{\gamma} \right)^{\left( \tfrac{\overline k_\gamma}{h} \right)} \max_{\nu \in \mc I}~w_\nu \|x^{0,\nu} - x^{*,\nu}\| \\
 & \geq \max_{\nu \in \mc I}~w_\nu \|x^{\overline k_\gamma,\nu} - x^{*,\nu}\|;
\end{align*}
\end{itemize}

\item In the proof of item (ii): the chain of inequalities reads as
\begin{align*}
 \displaystyle 
 & \max_{\nu \in \mc I}~w_\nu \|x^{\overline k_\gamma + 1,\nu} - x^{*,\nu}\| = \max_{\nu \in {\cal J}^{\overline k_\gamma}}~w_\nu \|x^{\overline k_\gamma + 1,\nu} - x^{*,\nu}\| \\
 & \qquad \leq \max_{\nu \in {\cal J}^{\overline k_\gamma}}~w_\nu \|T_\nu(\xbold^{\overline k_\gamma, -\nu}) - T_\nu(\xbold^{*, -\nu})\| + \\ 
 & \qquad \;\;\;\;\; \max_{\nu \in {\cal J}^{\overline k_\gamma}}~w_\nu \|\widetilde x^{\overline k_\gamma, \nu} - T_\nu(\xbold^{\overline k_\gamma, -\nu})\| + \\
 & \qquad \;\;\;\;\; \max_{\nu \in {\cal J}^{\overline k_\gamma}}~w_\nu \|T_\nu(\xbold^{*, -\nu}) - x^{*, \nu}\| +\\
 & \qquad \;\;\;\;\; \max_{\nu \in {\cal J}^{\overline k_\gamma}}~w_\nu \|x^{\overline k_\gamma + 1, \nu} - \widetilde x^{\overline k_\gamma, \nu}\| \\
 & \qquad \leq \max_{\nu \in {\cal J}^{\overline k_\gamma}}~w_\nu \|T_\nu(\xbold^{\overline k_\gamma, -\nu}) - T_\nu(\xbold^{*, -\nu})\| + \max_{\nu \in \mc I}~w_\nu (2 \beta \sqrt{i_\nu} + \varepsilon) \\
 & \qquad \leq \alpha \max_{\nu \in \mc I}~w_\nu \|x^{\overline k_\gamma,\nu} - x^{*,\nu}\| + \max_{\nu \in \mc I}~w_\nu (2 \beta \sqrt{i_\nu} + \varepsilon) \\
 & \qquad < \left( \tfrac{\gamma \alpha}{1 - \gamma \alpha} + 1 \right) \max_{\nu \in \mc I}~w_\nu (2 \beta \sqrt{i_\nu} + \varepsilon) = \left( \tfrac{1}{1 - \gamma \alpha} \right) \max_{\nu \in \mc I}~w_\nu (2 \beta \sqrt{i_\nu} + \varepsilon),
\end{align*}
where $\widetilde x^{\overline k_\gamma,\nu}$ is defined in \eqref{eq: inexact BR} for the iteration $\overline k_\gamma$;

\item The proof of item (iii) is akin to that of Theorem \ref{th: error bound}.
\end{itemize}

\section{Discussion on Assumptions \ref{as: block-contraction} and \ref{as: discrete bound}}\label{sec: assumptions}
We introduce now classes of \glspl{MI-NEP} that structurally meet the conditions in Assumptions \ref{as: block-contraction} and \ref{as: discrete bound}. For ease of reading, we will treat the two cases separately.

\subsection{Conditions for Assumption \ref{as: block-contraction} and their relations with strong monotonicity}\label{sec: assumption 1}

Let us assume that the cost functions $\theta_\nu$'s are $\mc C^2$ and define quantities: $\sigma_\nu \triangleq \inf_{\xbold \in X} \lambda_{\min} \left[ \nabla^2_{x^\nu x^\nu} \theta_\nu(\xbold) \right], \, \forall \nu \in \mc I$ and $\overline \sigma_{\nu\nu'} \triangleq \sup_{\xbold \in X} \left\| \nabla^2_{x^\nu x^{\nu'}} \theta_\nu (\xbold) \right\|_2, \, \forall  (\nu,\nu') \in \mc I^2$
where ${\lambda}_{\min}[A]$ is the smallest eigenvalue of the symmetric and positive semidefinite matrix $A$. In the spirit of \cite{scutari2012monotone}, we consider the ``condensed'' $N \times N$ real matrix $\Upsilon$, entry-wise defined as follows:
$$
\Upsilon_{\nu\nu'} \triangleq \left\{ \begin{array}{ll} \sigma_\nu, & \text{ if } \nu = \nu' \\ \overline \sigma_{\nu\nu'}, & \text{ otherwise} \end{array} \right., \qquad \forall \, (\nu,\nu') \in \mc I^2.
$$
Matrix $\Upsilon$ is strictly row diagonally dominant with weights $\wbold^{-1} \in \Re^{N}_{++}$ if, for all $\nu \in \mc I$, $w_\nu^{-1} \Upsilon_{\nu\nu} > \sum_{\nu' \in \mc I \setminus \{\nu\}} w_{\nu'}^{-1} \Upsilon_{\nu\nu'}$, i.e., $w_\nu^{-1} \sigma_\nu > \sum_{\nu' \in \mc I \setminus \{\nu\}} w_{\nu'}^{-1} \overline \sigma_{\nu\nu'}$.

\begin{proposition}\label{pr: sufficient for assumption 1}
 Let $\Upsilon$ be strictly row diagonally dominant with weights $\wbold^{-1} \in \Re^{N}_{++}$.
 Then, Assumption~\ref{as: block-contraction} holds true with weights $\wbold$ and modulus
 $$
 \alpha = \max_{\nu \in \mc I} \tfrac{\sum_{\nu' \in \mc I \setminus \{\nu\}} w_{\nu'}^{-1} \overline \sigma_{\nu\nu'}}{w_\nu^{-1} \sigma_\nu} < 1.
 $$
 \hfill$\square$
\end{proposition}
\textit{Proof.}
 For every $\zbold, \ybold \in X$ and $\nu \in \mc I$, from the optimality conditions we have:
 $$
 \begin{aligned}
 &\nabla_{x^\nu} \theta_\nu (T_\nu(\zbold^{-\nu}),\zbold^{-\nu})^{\scriptscriptstyle\top} (x^\nu - T_\nu(\zbold^{-\nu})) \geq 0, \quad \forall x^\nu \in X_\nu,\\
 &\nabla_{x^\nu} \theta_\nu (T_\nu(\ybold^{-\nu}),\ybold^{-\nu})^{\scriptscriptstyle\top} (x^\nu - T_\nu(\ybold^{-\nu})) \geq 0, \quad \forall x^\nu \in X_\nu.
  \end{aligned}
 $$
 Therefore, we have:
 $$
 \begin{aligned}
& \nabla_{x^\nu} \theta_\nu (T_\nu(\zbold^{-\nu}),\zbold^{-\nu})^{\scriptscriptstyle\top} (T_\nu(\ybold^{-\nu}) - T_\nu(\zbold^{-\nu})) \geq 0,\\
&\nabla_{x^\nu} \theta_\nu (T_\nu(\ybold^{-\nu}),\ybold^{-\nu})^{\scriptscriptstyle\top} (T_\nu(\zbold^{-\nu}) - T_\nu(\ybold^{-\nu})) \geq 0,
  \end{aligned}
 $$
 and then, for some $\xbold \in \left((T_\nu(\zbold^{-\nu}),\zbold^{-\nu}), (T_\nu(\ybold^{-\nu}),\ybold^{-\nu})\right)$: $0 \leq (T_\nu(\ybold^{-\nu}) - T_\nu(\zbold^{-\nu}))^{\scriptscriptstyle\top} \left( \nabla_{x^\nu} \theta_\nu (T_\nu(\zbold^{-\nu}),\zbold^{-\nu}) - \nabla_{x^\nu} \theta_\nu (T_\nu(\ybold^{-\nu}),\ybold^{-\nu}) \right) \overset{\text{(a)}}{=} (T_\nu(\ybold^{-\nu}) - T_\nu(\zbold^{-\nu}))^{\scriptscriptstyle\top} \nabla^2_{x^\nu x^\nu} \theta_\nu(\xbold) (T_\nu(\zbold^{-\nu}) - T_\nu(\ybold^{-\nu})) + (T_\nu(\ybold^{-\nu}) - T_\nu(\zbold^{-\nu}))^{\scriptscriptstyle\top} \left( \sum_{\nu' \neq \nu} \nabla^2_{x^\nu x^{\nu'}} \theta_\nu(\xbold) (z^{\nu'} - y^{\nu'}) \right)$,
 where (a) follows by the mean-value theorem. We thus obtain that: $w_\nu \sigma_\nu \|T_\nu(\ybold^{-\nu}) - T_\nu(\zbold^{-\nu})\|_2^2\leq w_\nu (T_\nu(\ybold^{-\nu}) - T_\nu(\zbold^{-\nu}))^{\scriptscriptstyle\top} \nabla^2_{x^\nu x^\nu} \theta_\nu(\xbold) (T_\nu(\ybold^{-\nu}) - T_\nu(\zbold^{-\nu}))\leq (T_\nu(\ybold^{-\nu}) - T_\nu(\zbold^{-\nu}))^{\scriptscriptstyle\top} \left( w_\nu \sum_{\nu' \neq \nu} \nabla^2_{x^\nu x^{\nu'}} \theta_\nu(\xbold) (z^{\nu'} - y^{\nu'}) \right)\leq \|T_\nu(\ybold^{-\nu}) - T_\nu(\zbold^{-\nu})\|_2 \left( w_\nu \sum_{\nu' \neq \nu} \overline \sigma_{\nu \nu'} \|z^{\nu'} - y^{\nu'}\|_2 \right)\leq \|T_\nu(\ybold^{-\nu}) - T_\nu(\zbold^{-\nu})\|_2 \left( w_\nu \sum_{\nu' \neq \nu} w_{\nu'}^{-1} \overline \sigma_{\nu \nu'} w_{\nu'} \|z^{\nu'} - y^{\nu'}\|_2 \right)\leq \|T_\nu(\ybold^{-\nu}) - T_\nu(\zbold^{-\nu})\|_2 \left( w_\nu \sum_{\nu' \neq \nu} w_{\nu'}^{-1} \overline \sigma_{\nu \nu'} \right) \max_{\nu' \in \mc I}~w_{\nu'} \|z^{\nu'} - y^{\nu'}\|_2$.
This finally yields the following inequality, which proves the statement:
 \begin{align*}
  w_{\nu} \|T_\nu(\ybold^{-\nu}) - T_\nu(\zbold^{-\nu})\|_2 \leq \max_{\lambda \in \mc I} \tfrac{\sum_{\nu' \in \mc I \setminus \{\lambda\}} w_{\nu'}^{-1} \overline \sigma_{\lambda\nu'}}{w_{\lambda}^{-1} \sigma_\lambda} \max_{\nu' \in \mc I} w_{\nu'} \|z^{\nu'} - y^{\nu'}\|_2. \tag*{}
 \end{align*}
%
\hfill$\blacksquare$

\noindent
Proposition \ref{pr: sufficient for assumption 1} shows that strict row diagonal dominance of matrix $\Upsilon$ with weights $\wbold^{-1} \in \Re^{N}_{++}$ is a sufficient condition for Assumption \ref{as: block-contraction} to hold true.
However, it is not immediate how to relate this fact with any monotonicity property for the \gls{MI-NEP}.
Let us define the standard game mapping $F:\R^n\to\R^n$ as follows:
$$
F(\xbold) \triangleq \left( \nabla_{x^\nu} \theta_\nu (\xbold) \right)_{\nu \in \mc I}.
$$
With a slight abuse of terminology, we say that the \gls{MI-NEP} is strongly monotone with constant $\mu > 0$ if
$F$ is strongly monotone with constant $\mu$, i.e.,
$$
(F(\xbold)-F(\ybold))^\top (\xbold - \ybold) \geq \mu \|\xbold - \ybold\|^2_2, \quad \forall \, \xbold, \ybold \in X.
$$
In general, there is not a direct relation between strong monotonicity and strict row diagonal dominance of $\Upsilon$, as illustrated in the following examples:

\begin{example}\label{ex: no diagonal dominance}
 Consider a \gls{MI-NEP} with 3 players characterized by scalar decision variables $x^1, x^2, x^3 \in \mathbb{R}$, and cost functions
 $\theta_1(x^1,x^2,x^3) = \tfrac{3}{2}(x^1)^2 + 2 x^1 x^2 + 2 x^1 x^3$, $\theta_2(x^1,x^2,x^3) = \tfrac{3}{2}(x^2)^2 + 2 x^2 x^1 + 2 x^2 x^3$, $\theta_3(x^1,x^2,x^3) = \tfrac{3}{2}(x^3)^2 + 2 x^3 x^1 + 2 x^3 x^2$.
 In this case, we have:
 $$
 F(x^1,x^2,x^3) = \begin{pmatrix} 3 & 2 & 2 \\ 2 & 3 & 2 \\ 2 & 2 & 3 \end{pmatrix} \begin{pmatrix} x^1 \\ x^2 \\ x^3 \end{pmatrix}, \qquad \Upsilon = \begin{pmatrix} 3 & 2 & 2 \\ 2 & 3 & 2 \\ 2 & 2 & 3 \end{pmatrix}.
 $$
 Therefore the \gls{MI-NEP} is strongly monotone with constant $\mu = 1$, but do not exist weights $\wbold^{-1} \in \Re^{3}_{++}$ such that $\Upsilon$ is strictly row diagonally dominant.
 \hfill$\square$
\end{example}

\begin{example}\label{ex: no strongly monotone}
 Consider a \gls{MI-NEP} involving 2 players with decision variables $x^1, x^2 \in \mathbb{R}$, and cost functions $\theta_1(x^1,x^2) = (x^1)^2 + x^1 x^2$, $\theta_2(x^1,x^2) = 5(x^2)^2 + 9 x^2 x^1$.
 In this case:
 $$
 F(x^1,x^2) = \begin{pmatrix} 2 & 1 \\ 9 & 10 \end{pmatrix} \begin{pmatrix} x^1 \\ x^2 \end{pmatrix}, \qquad \Upsilon = \begin{pmatrix} 2 & 1 \\ 9 & 10 \end{pmatrix}.
 $$
 Therefore, the matrix $\Upsilon$ is strictly row diagonally dominant with unitary weights, but the \gls{MI-NEP} is not monotone, in fact
 \begin{equation}
	  \tfrac{1}{2} \left( JF+JF^\top \right) = \begin{pmatrix} 2 & 5 \\ 5 & 10 \end{pmatrix} \not\succcurlyeq 0. \tag*{}
 \end{equation}
\hfill$\square$
\end{example}

\noindent
However, with the following result we show that it is always possible to suitably perturb a strongly monotone \gls{MI-NEP} in order to obtain strict row diagonal dominance of matrix $\Upsilon$, and then meet Assumption \ref{as: block-contraction}.

\begin{proposition}\label{pr: strong monotonicity}
 Let the \gls{MI-NEP} in \eqref{eq: prob} be strongly monotone with constant $\mu > 0$.
 Then, for any given $\overline \alpha < 1$ and $\overline \xbold \in X$, the \gls{MI-NEP} whose cost functions $\overline \theta_\nu$ for any $\nu \in \mc I$ are defined in one of the following ways:
 \begin{enumerate}[(i)]
  \item $\overline \theta_\nu (\xbold) \triangleq \theta_\nu (\xbold) + \tfrac{\eta_\nu}{2} \|x^\nu - \overline x^\nu\|_2^2$, with
  $
   \eta_\nu \triangleq \max \left\{ \tfrac{\sum_{\nu' \in \mc I \setminus \{\nu\}} \overline \sigma_{\nu\nu'}}{\overline \alpha} - \mu , 0 \right\},
  $
  \item $\overline \theta_\nu (\xbold) \!\triangleq\! \theta_\nu (\xbold) \!+\! \tfrac{\rho_\nu}{2} \|x^\nu\|^2_{\nabla^2_{x^\nu x^\nu} \theta_\nu (\overline \xbold)}$, with
  $
   \rho_\nu \!\triangleq\! \max \left\{ \tfrac{\sum_{\nu' \in \mc I \setminus \{\nu\}} \overline \sigma_{\nu\nu'}}{\overline \alpha \mu} \!-\! 1, 0 \right\},
  $
 \end{enumerate}
 is strongly monotone and meets Assumption \ref{as: block-contraction} with $\overline \alpha$ and $w_\nu = 1$ for all $\nu \in \mc I$.
 \hfill$\square$
\end{proposition}
\textit{Proof.}
Since $JF(\xbold) \succeq \mu I$ we have $\sigma_\nu \geq \mu$ for all $\nu \in \mc I$ and $\xbold \in X$.

Consider the condensed matrix $\overline \Upsilon$ related to the \gls{MI-NEP} with perturbed cost functions, and note that it differs \gls{wrt} $\Upsilon$ of the original \gls{MI-NEP} only in its diagonal elements.
 Specifically, if the $\nu$-th player problem is defined as in case (i), then the corresponding diagonal element is such that
 \begin{align*}
 \overline \Upsilon_{\nu \nu} &= \sigma_\nu + \eta_\nu \geq \mu + \eta_\nu \geq \mu + \tfrac{\sum_{\nu' \in \mc I \setminus \{\nu\}} \overline \sigma_{\nu\nu'}}{\overline \alpha} - \mu = \tfrac{\sum_{\nu' \in \mc I \setminus \{\nu\}} \overline \sigma_{\nu\nu'}}{\overline \alpha}.
 \end{align*}
 Otherwise, in case (ii), it holds that
 \begin{align*}
 \overline \Upsilon_{\nu \nu} &= \sigma_\nu + \rho_\nu \lambda_{\min} \left[ \nabla^2_{x^\nu x^\nu} \theta_\nu(\overline \xbold) \right] \geq (1 + \rho_\nu) \sigma_\nu \geq (1 + \rho_\nu) \mu \\
 &\geq \left(1 + \tfrac{\sum_{\nu' \in \mc I \setminus \{\nu\}} \overline \sigma_{\nu\nu'}}{\overline \alpha \mu} - 1 \right) \mu = \tfrac{\sum_{\nu' \in \mc I \setminus \{\nu\}} \overline \sigma_{\nu\nu'}}{\overline \alpha}.
 \end{align*}
 Therefore, in any case it holds that
 $
 \max_{\nu \in \mc I}~(\sum_{\nu' \in \mc I \setminus \{\nu\}} \overline \sigma_{\nu\nu'})/\overline \Upsilon_{\nu \nu} \leq \overline \alpha < 1,
 $
 and Proposition \ref{pr: sufficient for assumption 1} can be used to conclude the proof.
\hfill$\blacksquare$

\noindent
Thus, Proposition \ref{pr: strong monotonicity} shows two different ways to perturb any strongly monotone \gls{MI-NEP} to meet Assumption \ref{as: block-contraction} with any desired contraction constant $\overline \alpha$. 
Notice that in case (i) the perturbation considered is nothing else than a classical proximal term that is often employed in numerical methods. In a game-theoretic context, for instance, a similar result was already given for fully continuous problems \cite[Prop.~12.17]{FacchPang10}, while it has been extended here to a mixed-integer setting. On the other hand, the perturbation used in case (ii) introduces a quadratic term to strengthen the degree of strong monotonicity of the problem.
In case every $\eta_\nu$ or $\rho_\nu$ are equal to zero, then the original \gls{MI-NEP} clearly satisfies Assumption \ref{as: block-contraction} with $\overline \alpha$.

\subsection{On Assumption~\ref{as: discrete bound}}\label{sec: assumption 2}

Assuming the boundedness of each $\Omega_\nu$ implies the existence of some $\beta$ large enough so that Assumption~\ref{as: discrete bound} is met. However, the smaller the $\beta$, the tighter the bounds established in Theorems~\ref{th: error bound}, \ref{th: error bound continuous} and \ref{th: error bound wrt continuous NEP}. Therefore, an exceedingly large value for $\beta$ could yield irrelevant error bounds.
We thus introduce here some classes of \glspl{MI-NEP} for which Assumption~\ref{as: discrete bound} is met with a reasonably small $\beta$.

\begin{proposition}\label{pr: sufficient for beta general}
 \sloppy For all $\nu \in \mc I$ and $\xbold \in X$, suppose that $\nabla_{x^\nu} \theta_\nu (\cdot, \xbold^{-\nu})$ is Lipschitz continuous and strongly monotone with constants $L$ and $\sigma$, respectively. Moreover, assume that $X_\nu = [\overline l^\nu, \, \overline u^\nu] \times \widetilde X_\nu$, with $\overline l^\nu, \overline u^\nu \in \Z^{i_\nu}$ and $\widetilde X_\nu \subseteq \Re^{n_\nu-i_\nu}$.
 Then, Assumption \ref{as: discrete bound} is verified with $\beta = \tfrac{1}{2} \sqrt{L/\sigma}$ and the Euclidean norm.
 \hfill$\square$
\end{proposition}
\textit{Proof.}
 Let $\overline x^\nu$ be the feasible point defined by $\overline x^\nu_j=\lceil T_\nu(\xbold^{-\nu})_j \rceil$ if $\lceil T_\nu(\xbold^{-\nu})_j \rceil - T_\nu(\xbold^{-\nu})_j \leq T_\nu(\xbold^{-\nu})_j - \lfloor T_\nu(\xbold^{-\nu})_j \rfloor$, $\overline x^\nu_j=\lfloor T_\nu(\xbold^{-\nu})_j \rfloor$ otherwise,
 for any $j \in \{1,\ldots,i_\nu\}$, and $\overline x^\nu_j = T_\nu(\xbold^{-\nu})_j$ for any $j \in \{i_\nu + 1,\ldots,n_\nu\}$.
 Then, we have $\theta_\nu (\overline x^\nu, \xbold^{-\nu}) - \theta_\nu (T_\nu(\xbold^{-\nu}), \xbold^{-\nu})\leq \nabla_{x^\nu} \theta_\nu (T_\nu(\xbold^{-\nu}), \xbold^{-\nu})^{\scriptscriptstyle\top} (\overline x^\nu - T_\nu(\xbold^{-\nu})) + \tfrac{L}{2} \|\overline x^\nu - T_\nu(\xbold^{-\nu})\|_2^2= \tfrac{L}{2} \|\overline x^\nu - T_\nu(\xbold^{-\nu})\|_2^2 \leq (L/8) i_\nu$,
 where the first inequality follows by the descent lemma \cite[Prop.~A.24]{bertsekas1997nonlinear}, while the second equality is true since, if $|\overline x^\nu_j - T_\nu(\xbold^{-\nu})_j| \neq 0$, then 
$\nabla_{x^\nu_j} \theta_\nu (T_\nu(\xbold^{-\nu}), \xbold^{-\nu}) = 0$.
 It is now clear that any $\widehat x^\nu \in R_\nu(\xbold^{-\nu})$ shall satisfy
 $
 \theta_\nu (\widehat x^\nu, \xbold^{-\nu}) - \theta_\nu (T_\nu(\xbold^{-\nu}), \xbold^{-\nu}) \leq \theta_\nu (\overline x^\nu, \xbold^{-\nu}) - \theta_\nu (T_\nu(\xbold^{-\nu}), \xbold^{-\nu}) \leq (L/8) i_\nu.
 $
 Thus, the following chain of inequalities holds true: $\tfrac{\sigma}{2} \|\widehat x^\nu - T_\nu(\xbold^{-\nu})\|_2^2 \leq \nabla_{x^\nu} \theta_\nu (T_\nu(\xbold^{-\nu}), \xbold^{-\nu})^{\scriptscriptstyle\top} (\widehat x^\nu - T_\nu(\xbold^{-\nu})) + \tfrac{\sigma}{2} \|\widehat x^\nu - T_\nu(\xbold^{-\nu})\|_2^2 \leq \theta_\nu (\widehat x^\nu, \xbold^{-\nu}) - \theta_\nu (T_\nu(\xbold^{-\nu}), \xbold^{-\nu}) \leq (L/8) i_\nu$
 where the first inequality holds in view of the first order optimality condition, while the second one follows from the strong monotonicity of $\nabla_{x^\nu} \theta_\nu (\cdot, \xbold^{-\nu})$ -- see, e.g., \cite[\S B.1.1]{bertsekas2015convex}.
 Thus, we obtain that
 $
 \|\widehat x^\nu - T_\nu(\xbold^{-\nu})\|_2 \leq \tfrac{1}{2} \sqrt{(L~i_\nu/\sigma)},
 $
 which concludes the proof.
\hfill$\blacksquare$

\noindent
Let us now consider, instead, the following case, which yields a tighter bound.

\begin{proposition}\label{pr: sufficient for beta 1}
 For all $\nu \in \mc I$, let:
 $
 \theta_\nu (x^\nu, \xbold^{-\nu}) = \sum_{j = 1}^{i_\nu} \overline \theta_{\nu,j} (x^{\nu}_j, \xbold^{-\nu}) + \widetilde \theta_{\nu} \left( \left(x^{\nu}_j\right)_{j = i_\nu+1}^{n_\nu}, \xbold^{-\nu} \right), \; X_\nu = [\overline l^\nu, \, \overline u^\nu] \times \widetilde X_\nu,
 $
 where any $\overline \theta_{\nu,j} : \Re^{1+n-n_{\nu}} \to \Re$, $\widetilde \theta_{\nu} : \Re^{n-i_{\nu}} \to \Re$, $\overline l^\nu, \overline u^\nu \in \Re^{i_\nu}$, and $\widetilde X_\nu \subseteq \Re^{n_\nu-i_\nu}$.
 Then, Assumption~\ref{as: discrete bound} is verified with $\beta = 1$ and the Euclidean norm.
 \hfill$\square$
\end{proposition}
\textit{Proof.}
 In view of the convexity of any function $\overline \theta_{\nu,j} (\cdot, \xbold^{-\nu})$, we can conclude that, for any $\nu \in \mc I$ and any $\widehat x^\nu \in R_\nu(\xbold^{-\nu})$:
 $$
 |\widehat x^\nu_j - T_\nu(\xbold^{-\nu})_j| \leq \max \{ \lceil T_\nu(\xbold^{-\nu})_j \rceil - T_\nu(\xbold^{-\nu})_j, T_\nu(\xbold^{-\nu})_j - \lfloor T_\nu(\xbold^{-\nu})_j \rfloor\} \leq 1,
 $$
 for all $j \in \{1,\ldots,i_\nu\}$, while
 $
 \| \left(\widehat x^{\nu}_j\right)_{j = i_\nu+1}^{n_\nu} - \left(T_\nu(\xbold^{-\nu})_j\right)_{j = i_\nu+1}^{n_\nu} \|_2 = 0.
 $
 Then, we obtain $\|\widehat x^\nu - T_\nu(\xbold^{-\nu})\|_2 \leq \sqrt{i_\nu}$, and the result follows immediately.
\hfill$\blacksquare$

 \noindent
 Next, we identify more restrictive conditions producing an even tighter bound:

\begin{proposition}\label{pr: sufficient for beta 2}
 Assume the same setting as in Proposition \ref{pr: sufficient for beta 1}. For all $\nu \in \mc I$, suppose further that $\overline l^\nu, \overline u^\nu \in \Z^{i_\nu}$, and
 $
 \overline \theta_{\nu,j} (x^{\nu}_j, \xbold^{-\nu}) = \tfrac{1}{2} q^\nu_j(\xbold^{-\nu}) (x^{\nu}_j)^2 + c^\nu_j(\xbold^{-\nu})x^{\nu}_j, \forall \, j \in \{1,\ldots,i_\nu\},
 $
 where $q^\nu_j : \Re^{n-n_{\nu}} \to \Re_{++}$ and $c^\nu_j : \Re^{n-n_{\nu}} \to \Re$.
 Then, Assumption~\ref{as: discrete bound} is verified with $\beta = 1/2$ and the Euclidean norm.
 \hfill$\square$
\end{proposition}
\textit{Proof.}
 By exploiting the proof of Proposition \ref{pr: sufficient for beta 1},
 we only need to show that, for any $\nu \in \mc I$ and any $\widehat x^\nu \in R_\nu(\xbold^{-\nu})$, the following holds true:
 \begin{equation}\label{eq: bound 2 quad}
 |\widehat x^\nu_j - T_\nu(\xbold^{-\nu})_j| = \min \{ \lceil T_\nu(\xbold^{-\nu})_j \rceil - T_\nu(\xbold^{-\nu})_j, T_\nu(\xbold^{-\nu})_j - \lfloor T_\nu(\xbold^{-\nu})_j \rfloor\},
 \end{equation}
 for all $j \in \{1,\ldots,i_\nu\}$, since in this case we obtain $|\widehat x^\nu_j - T_\nu(\xbold^{-\nu})_j| \leq 1/2$.
 First, we observe that both $\lceil T_\nu(\xbold^{-\nu})_j \rceil$ and $\lfloor T_\nu(\xbold^{-\nu})_j \rfloor$ are in $[l^\nu_j,u^\nu_j]$. Moreover, if $T_\nu(\xbold^{-\nu})_j \in \Z$, then any $\widehat x^\nu_j$ must be equal to $T_\nu(\xbold^{-\nu})_j$, since $q^\nu_j(\xbold^{-\nu}) > 0$. 
 Therefore, we only have to consider the case $T_\nu(\xbold^{-\nu})_j \notin \Z$, which implies
 $
 q^\nu_j(\xbold^{-\nu}) (T_\nu(\xbold^{-\nu})_j) + c^\nu_j(\xbold^{-\nu}) = 0,
 $
 since $u^\nu_j, l^\nu_j \in \Z$. We thus have that $ \tfrac{1}{2} q^\nu_j(\xbold^{-\nu}) (x^{\nu}_j)^2 + c^\nu_j(\xbold^{-\nu})x^{\nu}_j = \tfrac{1}{2} q^\nu_j(\xbold^{-\nu}) (T_\nu(\xbold^{-\nu})_j)^2 + c^\nu_j(\xbold^{-\nu})T_\nu(\xbold^{-\nu})_j + \tfrac{1}{2} q^\nu_j(\xbold^{-\nu}) (x^{\nu}_j - T_\nu(\xbold^{-\nu})_j)^2$.
Since $\widehat x^\nu_j$ is an integer minimizer of this univariate quadratic function, it shall be the closest to $T_\nu(\xbold^{-\nu})_j$ because $q^\nu_j(\xbold^{-\nu}) > 0$.
 This implies \eqref{eq: bound 2 quad}, and hence we obtain $\|\widehat x^\nu - T_\nu(\xbold^{-\nu})\|_2 \leq (1/2) \sqrt{i_\nu}$.
\hfill$\blacksquare$

\noindent
In Propositions~\ref{pr: sufficient for beta general}--\ref{pr: sufficient for beta 2}, the continuous set $X_\nu = [\overline l^\nu, \, \overline u^\nu] \times \widetilde X_\nu$ of each player has a separable structure. Removing this condition is not reasonable to meet Assumption~\ref{as: discrete bound} with a suitable bound $\beta$, as supported by the following example:

\begin{example}\label{ex: separable feasible set}
 \sloppy Let $N=1$ with $\theta_1(x^1) = (x^1_1)^2 + (x^1_2)^2$, $i_1 = 2$, $X_1 = \left\{x^1 \in \Re^2 \mid x^1_1 \geq \tfrac{1}{2}, \, x^1_2 \geq \upsilon x^1_1 - \tfrac{\upsilon-1}{2} \right\}$. By considering any $\upsilon \geq 1$, it holds that $T_1 = \left(\tfrac{1}{2},\tfrac{1}{2}\right)^{\scriptscriptstyle\top}$, and $\widehat x^1 = \left(1, \tfrac{\upsilon+1}{2} \right)^{\scriptscriptstyle\top}$. Therefore, the distance $\|\widehat x^1 - T_1\|_2 = \tfrac{\sqrt{\upsilon^2+1}}{2}$ appearing in Assumption~\ref{as: discrete bound} depends on $\upsilon$ and can be arbitrarily large.
 \hfill$\square$
\end{example}

\section{Practical usage of the BR algorithms and numerical results}\label{sec: numerical}

The results developed in this paper allow one to make use (or combine) both Algorithms \ref{alg: Jacobi} and \ref{alg: Jacobi continuous} for the computation of \gls{MI-NE} for the class of \glspl{MI-NEP} satisfying Assumptions \ref{as: block-contraction} and \ref{as: discrete bound}. Specifically, consider the following procedures:

\begin{enumerate}[(i)]
\item Using Algorithm \ref{alg: Jacobi} only. This procedure does not have any theoretical guarantee of success, however, if convergence happens, then it certainly returns a solution. In any case, combining Theorems \ref{th: error bound} and \ref{th: error bound wrt continuous NEP} allows us to conclude that the sequence $\{\xbold^k\}_{k\in\mathbb{N}}$: i) belongs to a region, whose diameter depends on $\alpha$ and $\beta$, that contains every possible solution, and ii) it is bounded, even if $\Omega$ is unbounded and the \gls{MI-NEP} does not admit any solution.

\item Using Algorithm \ref{alg: Jacobi continuous} to compute the unique solution of the relaxed \gls{NEP}, $\overline{\xbold}$, and then, starting from this point, using Algorithm \ref{alg: Jacobi} to compute a \gls{MI-NE}. By Theorem \ref{th: error bound continuous}, $\overline{\xbold}$ shall be reasonably close to the solution set $\mathcal S$ of the \gls{MI-NEP} (according to the values of $\alpha$, $\beta$), and it can be computed almost inexpensively -- see Theorem \ref{th: error bound continuous in continuous problems}. The considerations in (i) also apply here, however Algorithm \ref{alg: Jacobi} could benefit from starting closer to $\mathcal S$.

\item Using Algorithm \ref{alg: Jacobi} to compute a reduced feasible region around the solution set $\mathcal{S}$ of the \gls{MI-NEP}, and then using an enumerative method over such a reduced region (see Section \ref{sec: introduction} for references) to compute a solution (or $\mathcal{S}$ itself). 
According to Theorem \ref{th: error bound}, which gives theoretical guarantees of convergence for this procedure, the ratio between the size of the reduced region and the original feasible set depends on $\alpha$ and $\beta$, thus strongly affecting the performance of the enumerative method employed.

\item Using Algorithm \ref{alg: Jacobi continuous} to compute a reduced feasible region around the solution set $\mathcal{S}$ of the \gls{MI-NEP}, and then using an enumerative method over such a reduced region to compute a solution (or $\mathcal{S}$ itself). In addition to the same comments in (iii), note that Algorithm \ref{alg: Jacobi continuous} is in general more efficient than Algorithm \ref{alg: Jacobi}, and the bound produced with this procedure (Theorem \ref{th: error bound continuous}) is better than that produced in item (iii) (Theorem \ref{th: error bound}).
\end{enumerate}

\noindent
%
We next compare procedures (i) and (ii) above, while enumerative methods as in (iii) and (iv) will be analyzed in future works. Specifically, we verify our findings on a numerical instance of a smart building control application.

\subsection{Problem description: Local smart building control}
Inspired by game-theoretic approaches to smart grids control applications \cite{cenedese2019charging}, we consider a smart building consisting of $N$ units (i.e., users, indexed by the set $\mc{I} \triangleq \{1, \ldots, N\}$) where each one of them is interested in designing an optimal schedule to switch on/off $m_\nu$ high power domestic appliances $\mc{A}_\nu \triangleq \{1, \ldots, m_\nu\}$ (e.g., washing machines, dishwashers, tumble dryers, electric vehicles) with known amount of required energy $\bar{u}^{h}_\nu > 0$, $h \in \mc{A}_\nu$ and $\nu \in \mc{I}$, over some time window $\mc{T} \triangleq \{1, \ldots, T\}$ to make the energy supply of the building smart and efficient. To this end, we assume each $\nu \in \mc{I}$ endowed with some storing capacity as, e.g., a battery or the electric vehicle itself.

The scheduling decision variable of each user consists of an integer vector $\delta^\nu \in \Delta_\nu^{T m_\nu}$ denoting the percentage of utilization of a certain appliance in period $k \in \mc{T}$, while the continuous one $u^\nu \in [0, u^\textrm{max}_\nu]^T$ regulates the acquisition of energy over $\mc{T}$, with $u^\textrm{max}_\nu > 0$. We then consider a scenario in which each single user has an individual supply contract with cost per unit $p_\nu > 0$ over the whole of $\mc{T}$. The local cost incurred by user $\nu \in \mc{I}$ can be formalized as:
\begin{equation}\label{eq:local_cost_function_mod}
    \begin{aligned}
    &\theta_\nu(z^\nu, \zbold^{-\nu}) =\\
    &\sum_{k \in \mc{T}}^{} \left[\kappa_\nu u^\nu(k)^2 \!+\! \chi_\nu \delta^\nu(k)^2 \!+\! p_\nu(k) \, a(\ubold(k)) u^\nu(k) \!+\! c_\nu \sum_{h \in \mc{A}_\nu} \left(y^\nu_h(k) \!-\! \delta_h^\nu(k) \bar{u}^{h}_\nu\right)^2 \right]\,,
    \end{aligned}
\end{equation}
where $\kappa_\nu > 0$ penalizes unnecessary energy acquisitions from the grid through the quadratic term ${u^\nu(k)}^2$, $\chi_\nu > 0$ favours low powers cycles through the quadratic term ${\delta^\nu(k)}^2$, and $p_\nu(k)$ reflects possible different tariffs across the day (daily vs night price of energy), while $a(\ubold(k)) \triangleq \sum_{\nu \in \mc{I}}^{} u^\nu(k)$ denotes the aggregate demand of energy associated with the set of users at time $k$, $z^\nu \triangleq \textrm{col}(u^\nu, \delta^\nu, y^\nu)$. Finally, the parameter $c_\nu > 0$ penalizes the deviation of some continuous variable $y^\nu_h \in \R_{+}^{T m_\nu}$ from the actual energy consumption for switching on a certain appliance given by $\delta_h^\nu(k) \bar{u}^{h}_\nu$. This last term represents a soft constraint forcing equality $y^\nu_h(k) = \delta_h^\nu(k) \bar{u}^{h}_\nu$, and hence the auxiliary variable $y^\nu_h$ acts as a proxy for the consumption required by appliance $h \in \mc{A}_\nu$.

Both the scheduling and energy acquisition variables are also subject to private constraints. For instance, we may assume that each appliance has to complete its task over the whole of $\mc{T}$, and this translates into:
$
	\sum_{k \in \mc{T}} y_h^\nu(k) = \bar\Delta_\nu \bar{u}^{h}_\nu, \text{ for all } h \in \mc{A}_\nu \,,
$
where $\bar\Delta_\nu$ is the largest element in $\Delta_\nu$, for all $\nu \in \mc I$.
Then, if $x_\nu(0) \ge 0$ denotes the initial \gls{SOC} of each storage unit, one has to satisfy for all $k \in \mc{T}$, $x_\nu(k+1) = x_\nu(k) + \eta_\nu u^\nu(k) - (\xi_\nu/\bar\Delta_\nu) \sum_{h \in \mc{A}_\nu} y^\nu_h(k)$ and $x_\nu(k) \geq 0$,
where $\eta_\nu$, $\xi_\nu > 0$ are some positive parameters representing the charging/discharging efficiency. 
To account for a possible physical cap $\bar{u}^\textrm{max}_\nu > 0$ limiting the delivery of energy in each time period, we shall also impose that 
$
	\sum_{h \in \mc{A}_\nu} y^\nu_h(k) \le \bar\Delta_\nu \bar{u}^\textrm{max}_\nu, \, \text{ for all } k \in \mc{T} \, .
$
With $\mc{Z}_\nu \triangleq [0, u^\textrm{max}_\nu]^T \times \Delta_\nu^{T m_\nu} \times \R_{+}^{T m_\nu}$, the resulting \gls{MI-NEP} thus reads as:
\begin{equation}\label{eq:mi_nep}
	\forall \nu \in \mc{I} : \left\{
	\begin{aligned}
		& \underset{z^\nu \in \mc{Z}_\nu}{\textrm{min}} && \theta_\nu(z^\nu,\zbold^{-\nu})\\
		&~\textrm{ s.t. } && \textstyle\sum_{k \in \mc{T}} y_h^\nu(k) = \bar\Delta_\nu \bar{u}^{h}_\nu, \forall h \in \mc{A}_\nu \,,\\
        &&& x_\nu(k\!+\!1) \!=\! x_\nu(k) \!+\! \eta_\nu u^\nu(k) \!-\! \tfrac{\xi_\nu}{\bar\Delta_\nu} \textstyle\sum_{h \in \mc{A}_\nu} y^\nu_h(k) \, \forall k \in \mc{T},\\
        &&& x_\nu(k) \geq 0, \; \textstyle\sum_{h \in \mc{A}_\nu} y^\nu_h(k) \le \bar\Delta_\nu \bar{u}^\textrm{max}_\nu, \, \forall k \in \mc{T}
	\end{aligned}
	\right.
\end{equation}
The final \gls{MI-NEP} turns out to be quadratic with asymmetries due to the different energy prices $p_\nu$ across agents. According to Proposition~\ref{pr: strong monotonicity}, we note that to make the \gls{MI-NEP} diagonally-dominant it suffices to adjust (specifically, increase) the design parameters $\kappa_\nu$, $\chi_\nu$ and/or $c_\nu$, since this would have the same effect on the considered costs in \eqref{eq:local_cost_function_mod} to having a proximal-like term without reference (i.e., $\overline x^\nu = 0$ in Proposition~\ref{pr: strong monotonicity}.(i)). The same consideration also applies to the case considered in Proposition~\ref{pr: strong monotonicity}.(ii), where one would simply need to know some feasible $\overline{\zbold}$ associated to the relaxed problem to compute $\nabla^2_{z^\nu z^\nu} \theta_\nu (\overline \zbold)$, which on the other hand also depends on some chosen values for $\kappa_\nu$, $\chi_\nu$ and $c_\nu$ themselves. In both cases, one can thus define a-priori some desired level of contraction $\overline \alpha <1$ to form a basis for the design of the weights appearing in the costs \eqref{eq:local_cost_function_mod}. For a careful choice of the latter parameters, note that one would ideally recover the existence (and uniqueness) condition established in Proposition~\ref{pr: existence}. In the limit case of an exceedingly low $\overline \alpha$, i.e., for too large values for $\kappa_\nu$, $\chi_\nu$ and $c_\nu$, however, one would obtain an almost decoupled problem (as the price $p_\nu$ is typically fixed and can not be manipulated) which inevitably could become of little significance from an application perspective. We will investigate these tight relations in future works.

\subsection{Numerical results}

All the experiments are carried out in Matlab on a laptop with an Apple M2 chip featuring 8-core CPU and 16 Gb RAM. The code has been developed in YALMIP environment \cite{Lofberg2004} with Gurobi \cite{gurobi} as solver to handle \glspl{MINLP}. 

We consider an instance of the \gls{MI-NEP} described in \eqref{eq:mi_nep} with $N = 8$ users willing to obviate the energy procurement over an horizon $T = 6$ and compute an optimal schedule to switch on/off $m_\nu \sim \mc{U}(2,4) \cap \Z$ domestic appliances ($\mc{U}(a,b)$ denotes the uniform distribution on the interval $[a,b]$). 
With $u_\nu^\textrm{max} = u^\textrm{max} = 1.2$, $x_\nu(0) = 0$, $c_\nu \sim 10^2 \cdot \mc{N}(6,0.02)$ (where $\mc{N}(a,b)$ denotes the normal distribution with mean $a$ and variance $b$), $\chi_\nu=c_\nu/10^2$, $\bar{u}^h_\nu \sim \mc{U}(1.2,2)$ and $\bar{u}^\textrm{max}_\nu = 1.2 \cdot \textrm{max}_h \, \{\bar{u}^h_\nu\}$, the individual price of energy follows two normal distributions to reflect daily and night tariffs, i.e., $\mc{N}((0.5+0.1)/(N u^\textrm{max}), 10^{-2})$ and $\mc{N}((0.35+0.05)/(N u^\textrm{max}), 10^{-2})$, respectively, while the parameter $\kappa_\nu \sim \mc{N}(6, 0.2)$. Specifically, we split the horizon $T$ in two parts: three hours associated to the daily consumption, and three to the night one.

According to the granularity specified for the set $\Delta_\nu$, we then conduct several numerical experiments. In particular, we will generically refer to \gls{MI-NEP}$^\text{u}$ to the case $\Delta_\nu = \{0,1,\ldots,100\}^{T m_\nu}$, i.e., the scheduling variable $\delta^\nu$ is allowed to take any integer value between $0$ and $100$, whereas we will refer to \gls{MI-NEP}$^\text{t}$ when $\Delta_\nu = \{0,10,\ldots,100\}^{T m_\nu}$, namely $\delta^\nu$ can only assume values corresponding to the tens between $0$ and $100$. These two cases have a different impact on the bound in \eqref{eq: error bound lim continuous}. From a random numerical instance of the considered \gls{MI-NEP}, for example, in view of the structure of the cost function in \eqref{eq:local_cost_function_mod}, we have $\alpha = 0.02$, $\beta = \tfrac{1}{2} \sqrt{L/\sigma} = 5.01$ (according to Proposition~\ref{pr: sufficient for beta general}), where $L = 1.2\cdot10^3$ and $\sigma = 12$ denote the Lipschitz constant of the (affine) game mapping and associated constant of strong monotonicity, respectively, which have been computed starting from the linear term characterizing the game mapping itself. Setting $w_\nu = w = 1$ for all $\nu \in \mc{I}$, the resulting bound in the RHS of \eqref{eq: error bound lim continuous} is $25.43$ for \gls{MI-NEP}$^\text{u}$, while it coincides with the same value multiplied by $100$ for \gls{MI-NEP}$^\text{t}$. For the application considered, while on the one hand taking a granularity up to the units may be restrictive from a practical perspective, the bound in the RHS of \eqref{eq: error bound lim continuous} is relevant to speed up the computation of a \gls{MI-NE}, whereas for the case \gls{MI-NEP}$^\text{t}$, albeit more realistic, the obtained bound is not meaningful. In this latter case, we therefore limit to propose a possible heuristic for improving the computation of an associated \gls{MI-NE}, for which however we do not have firm theoretical guarantees in the spirit of Theorem~\ref{th: error bound continuous}. 
With this regard, we will hence identify with \gls{MI-NEP}$_\text{r}$ those problems referring to the reduced feasible set, and making use of $\overline{x}^\nu$ as starting point for our algorithms.
On the contrary, \gls{MI-NEP}$_\text{f}$ will denote those examples considering the whole feasible set, initialized with $\xbold^0=0$. 

    Notice the abuse of notation in referring to those instances considering the reduced feasible set with $\overline{x}^\nu$ as starting point. Specifically, for the case with units we actually apply the bound in \eqref{eq: error bound lim continuous} around $\overline{x}^\nu$, thus running the algorithms onto a reduced feasible set. For the case with tens, instead, we heuristically observe a-posteriori the same behaviour as per the case with units, since we do not have any theoretical guarantee for the same bound.

For each randomly generated numerical instance of \eqref{eq:mi_nep}, we will thus end up exploring five different cases: \gls{MI-NEP}$^\triangle_\ast$, $\triangle \in \{\text{u}, \text{t}\}$, $\ast \in \{\text{f}, \text{r}\}$, plus the associated relaxed, continuous problem. According to Theorem~\ref{th: error bound continuous in continuous problems}, this latter always admits a unique Nash equilibrium that is computable via Algorithm~\ref{alg: Jacobi continuous}. We will finally contrast the performance of a Jacobi-type scheme ($\mc J^k=\mc I$ for all $k$, denoted as \texttt{J}) with a Gauss-Seidel one (players taken sequentially, one per iteration and denoted as \texttt{GS}). 

\begin{table}
	\caption{Computational time and number of iterations}
	\label{tab:num_res_J}
		\begin{tabular}{lllll}
			\toprule
			  & \gls{MI-NEP}$^\text{u}$$_\text{f}$  & \gls{MI-NEP}$^\text{u}$$_\text{r}$ & \gls{MI-NEP}$^\text{t}$$_\text{f}$ &  \gls{MI-NEP}$^\text{t}$$_\text{r}$\\
			\midrule
			\texttt{J}--CPU\textsubscript{eq}  & 4.66 [s] & 4.18 [s] & 5.87 [s] &  4.63 [s]\\
                \midrule
               \texttt{J}--\#Iter\textsubscript{eq} & 12.68 & 10.24 & 13.75 & 10.77\\
               \midrule
               \midrule
               \texttt{GS}--CPU\textsubscript{eq}  & 3.45 [s] & 3.30 [s] & 4.75 [s] & 3.70 [s]\\
               \midrule
               \texttt{GS}--\#Iter\textsubscript{eq} & 9.37 & 8.33 & 11.08 & 8.80\\
			\bottomrule
		\end{tabular}
\end{table}


We hence test our theoretical findings over $500$ random numerical instances of the \gls{MI-NEP} in \eqref{eq:mi_nep}. 
Specifically,
Table~\ref{tab:num_res_J} reports the average computational time and number of iterations needed for computing an \gls{MI-NE} in all those numerical examples in which Algorithm~\ref{alg: Jacobi} has converged, both for the Jacobi and Gauss-Seidel schemes.
For the experiments \gls{MI-NEP}$^\text{u}_\text{f}$, we have obtained a bound \eqref{eq: error bound lim continuous} always within the interval $[22.01, 25.61]$, with average value of $25.31$. 
As expected, running Algorithm~\ref{alg: Jacobi} to find a \gls{MI-NE} onto a reduced feasible set (columns \gls{MI-NEP}$^\triangle_\text{r}$, $\triangle \in \{\text{u}, \text{t}\}$) is faster than computing an equilibrium over the original feasible set (columns \gls{MI-NEP}$^\triangle_\text{f}$, $\triangle \in \{\text{u}, \text{t}\}$). In particular, while one can save around the $13$\% of the computational time in \gls{MI-NEP}$^\text{u}_\ast$, $\ast \in \{\text{f}, \text{r}\}$, since the reduction procedure brings the total number of integer variables approximately from $100^{144}$ to $50^{144}$ on average (as each $m_\nu \sim \mc{U}(2,4) \cap \Z$), this percentage grows markedly when considering the heuristic for the cases \gls{MI-NEP}$^\text{t}_\ast$, $\ast \in \{\text{f}, \text{r}\}$. 
From our numerical experience, computing the Nash equilibrium for the relaxed \gls{NEP}, which is always possible via Algorithm~\ref{alg: Jacobi continuous} in view of Theorem~\ref{th: error bound continuous in continuous problems}, takes $5.98$ iterations on average and it is extremely fast.

\begin{table}
	\caption{Failures}
	\label{tab:failures_J}
		\begin{tabular}{lllll}
			\toprule
			  & \gls{MI-NEP}$^\text{u}$$_\text{f}$  & \gls{MI-NEP}$^\text{u}$$_\text{r}$ & \gls{MI-NEP}$^\text{t}$$_\text{f}$ & \gls{MI-NEP}$^\text{t}$$_\text{r}$\\
			\midrule
			\texttt{J}--\% of failure  & 59.14 & 51.36 & 58.95 & 50.00\\
			\midrule
                \midrule
			\texttt{GS}--\% of failure  & 50.95 & 39.12 & 40.27 & 31.78\\
			\bottomrule
		\end{tabular}
\end{table}


\begin{figure}
    \centering
    \includegraphics[width=.5\columnwidth]{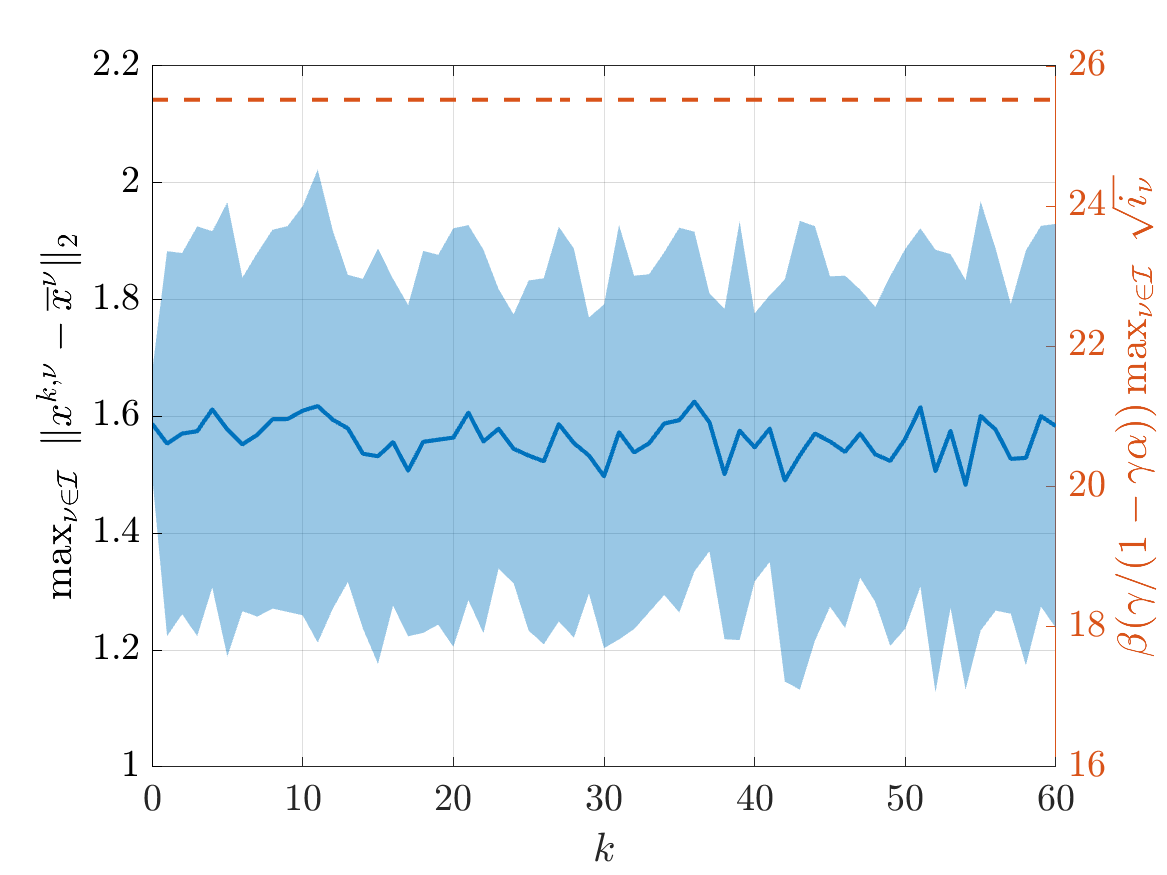}  
    \caption{Mean value (solid blue line) and standard deviation (shaded blue area) of the players' worst-case distance from the continuous equilibrium of the relaxed \gls{NEP}, computed in those numerical instances for which Algorithm~\ref{alg: Jacobi} has not converged for \gls{MI-NEP}$^\text{u}_\text{f}$. The red line denotes, instead, the bound in \eqref{eq: error bound wrt continuous NEP} with $\gamma=1.001$ (note the  different values on the left-right ordinates).}
    \label{fig:th4}
\end{figure}

While this analysis shows the practical impact of combining Algorithm~\ref{alg: Jacobi} and \ref{alg: Jacobi continuous} with the bound offered in Theorem~\ref{th: error bound continuous}, however, we should note that the \gls{BR} method as described in Algorithm~\ref{alg: Jacobi} may not necessarily produce a convergent sequence, neither in its Jacobi nor Gauss-Seidel versions (even though the \gls{MI-NEP} admits at least a \gls{MI-NE}). 
In particular, Table~\ref{tab:failures_J} shows the percentage of failures declared after $60$ iterations of both versions of Algorithm~\ref{alg: Jacobi} (actually, $60\cdot N=240$ for the Gauss-Seidel implementation) without computing an \gls{MI-NE}, established when two consecutive iterations meet the stopping criterion $\|\zbold^{k+1}-\zbold^{k}\|_2\le10^{-6}$. 
In general, we note that the percentage of failures associated to \gls{MI-NEP}$^\text{t}_\text{f}$ is lower than that of \gls{MI-NEP}$^\text{u}_\text{f}$. This can be explained as the case with units structurally admits way more combinations of integer feasible points compared to the one with tens, and this behaviour immediately reflects onto the cases considering reduced feasible sets, \gls{MI-NEP}$^\triangle_\text{r}$, $\triangle \in \{\text{u},\text{t}\}$. In addition, note that the high failure rate shown for \gls{MI-NEP}$^\text{u}_\ast$, $\ast \in \{\text{f}, \text{r}\}$ may be associated with the threshold employed to declare non-convergence, i.e., $60$ iterations. 
Despite computing a \gls{MI-NE} via Algorithm~\ref{alg: Jacobi} with a finer granularity is slightly faster, according to Table~\ref{tab:num_res_J}, it is also more likely to fail. 
Adopting Algorithm~\ref{alg: Jacobi} in combination with \ref{alg: Jacobi continuous} and the bound in Theorem~\ref{th: error bound continuous} may thus help in reducing the failures potentially occurring when Algorithm~\ref{alg: Jacobi} is used alone. Considering only those examples for \gls{MI-NEP}$^\text{u}_\text{f}$ in which convergence has not happened allows us to verify also the bound in Theorem~\ref{th: error bound wrt continuous NEP}, as shown in Fig.~\ref{fig:th4} for $\mc J^k=\mc I$ -- similar results are obtained for the Gauss-Seidel version.
\begin{figure}
    \centering
	\begin{subfigure}{.329\textwidth}
		\centering
		\includegraphics[width=\columnwidth]{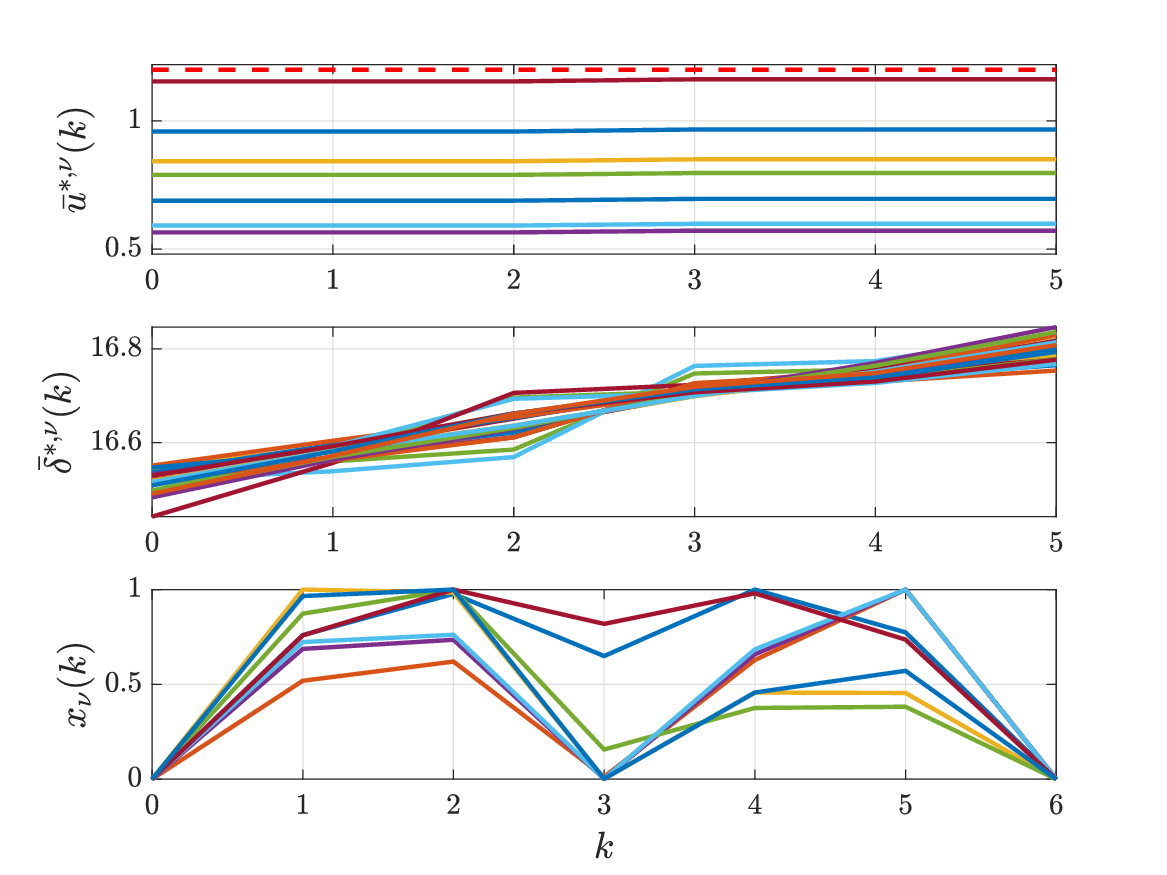}  
		\caption{}
		\label{fig:sub-first}
	\end{subfigure}~
	\begin{subfigure}{.329\textwidth}
		\centering
		\includegraphics[width=\columnwidth]{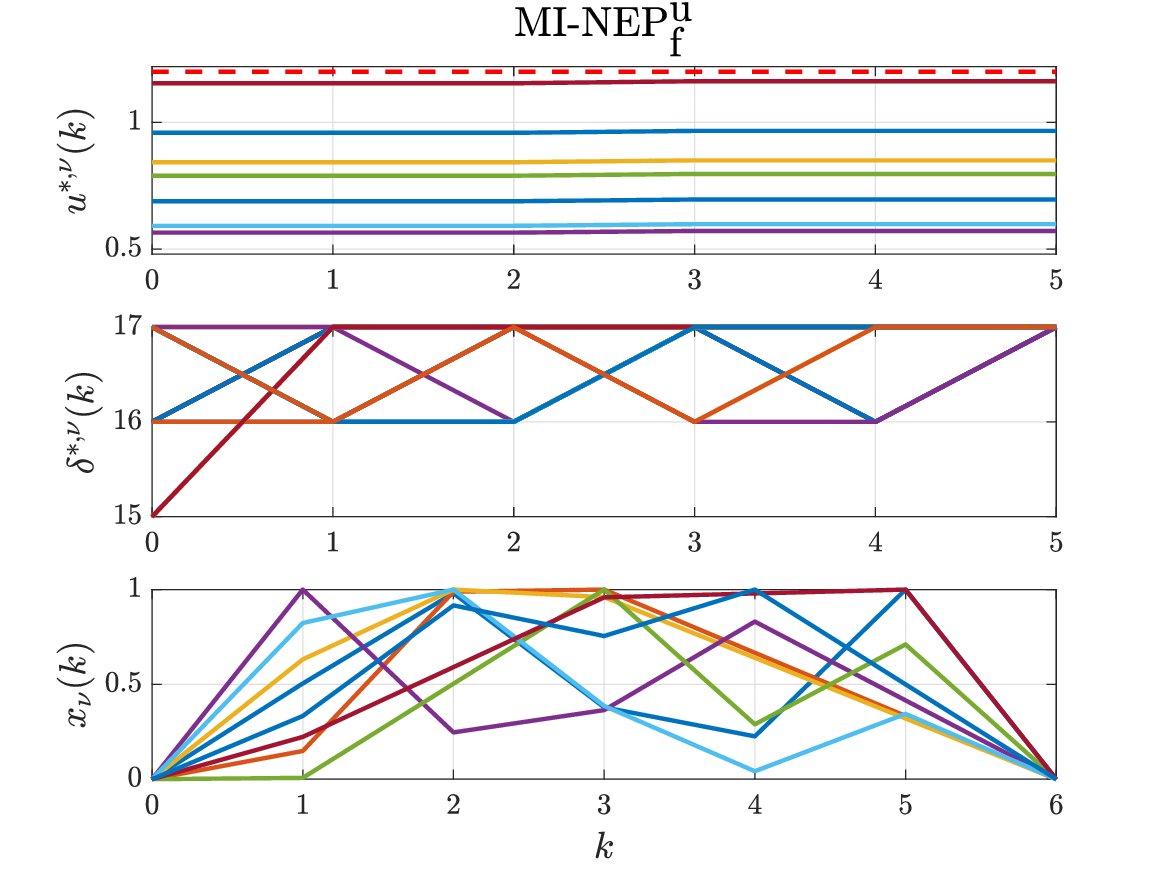}  
		\caption{}
		\label{fig:sub-second}
	\end{subfigure}~
        \begin{subfigure}{.329\textwidth}
		\includegraphics[width=\columnwidth]{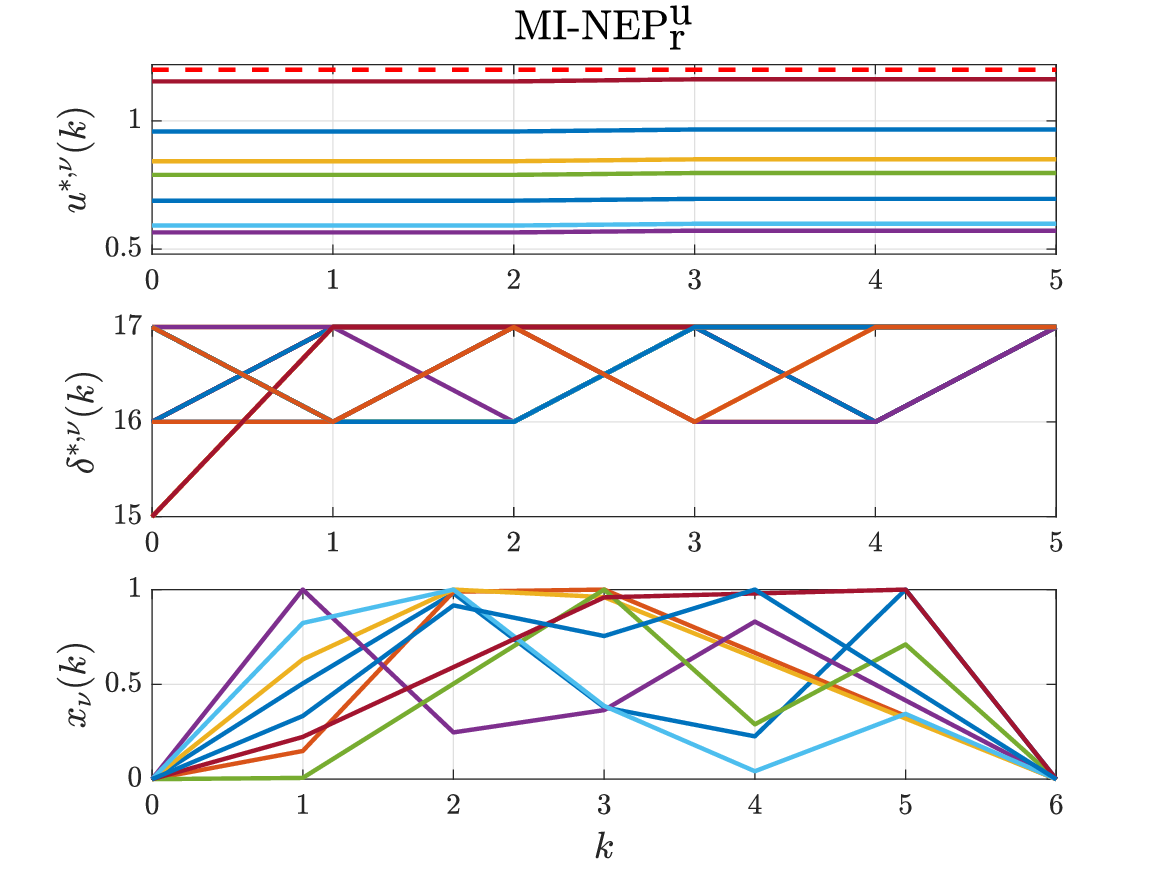}  
		\caption{}
		\label{fig:sub-third}
	\end{subfigure}
        \begin{subfigure}{.329\textwidth}
		\includegraphics[width=\columnwidth]{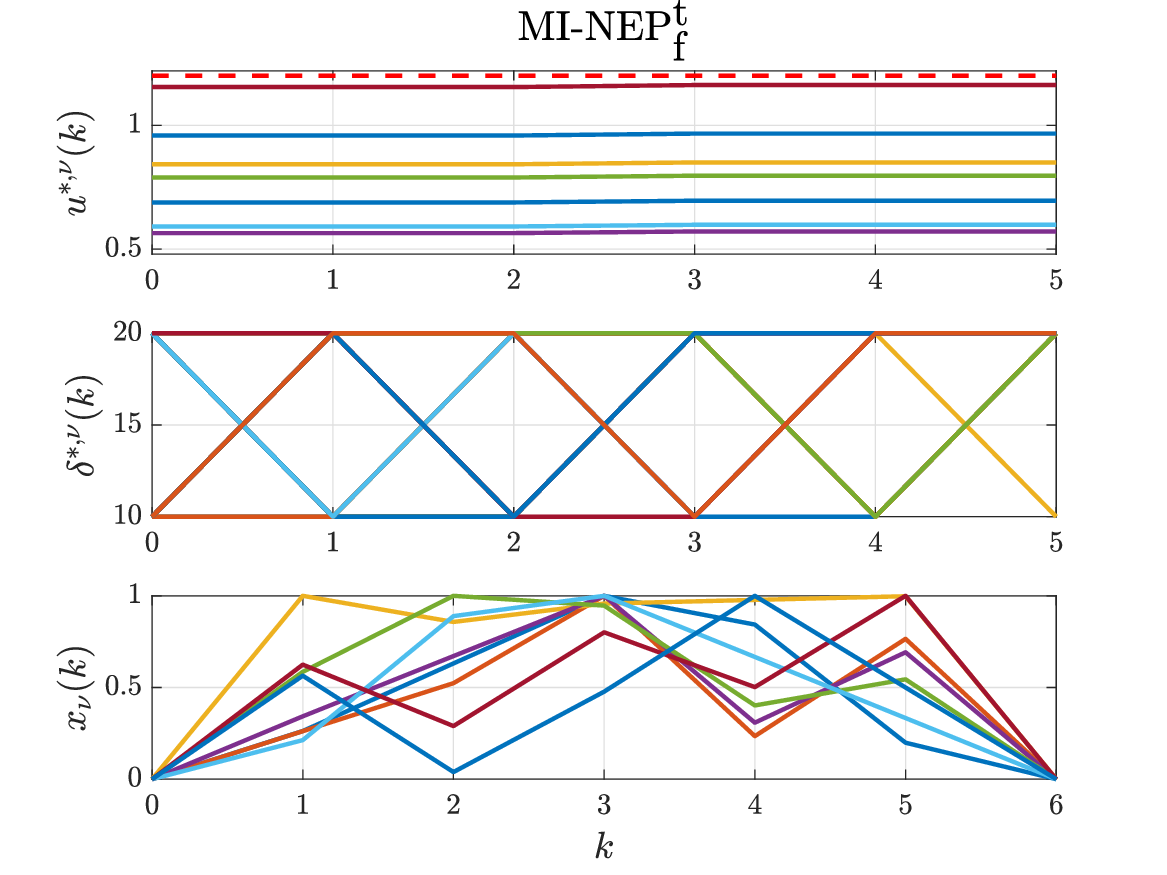}  
		\caption{}
		\label{fig:sub-fourth}
	\end{subfigure}~
        \begin{subfigure}{.329\textwidth}
		\includegraphics[width=\columnwidth]{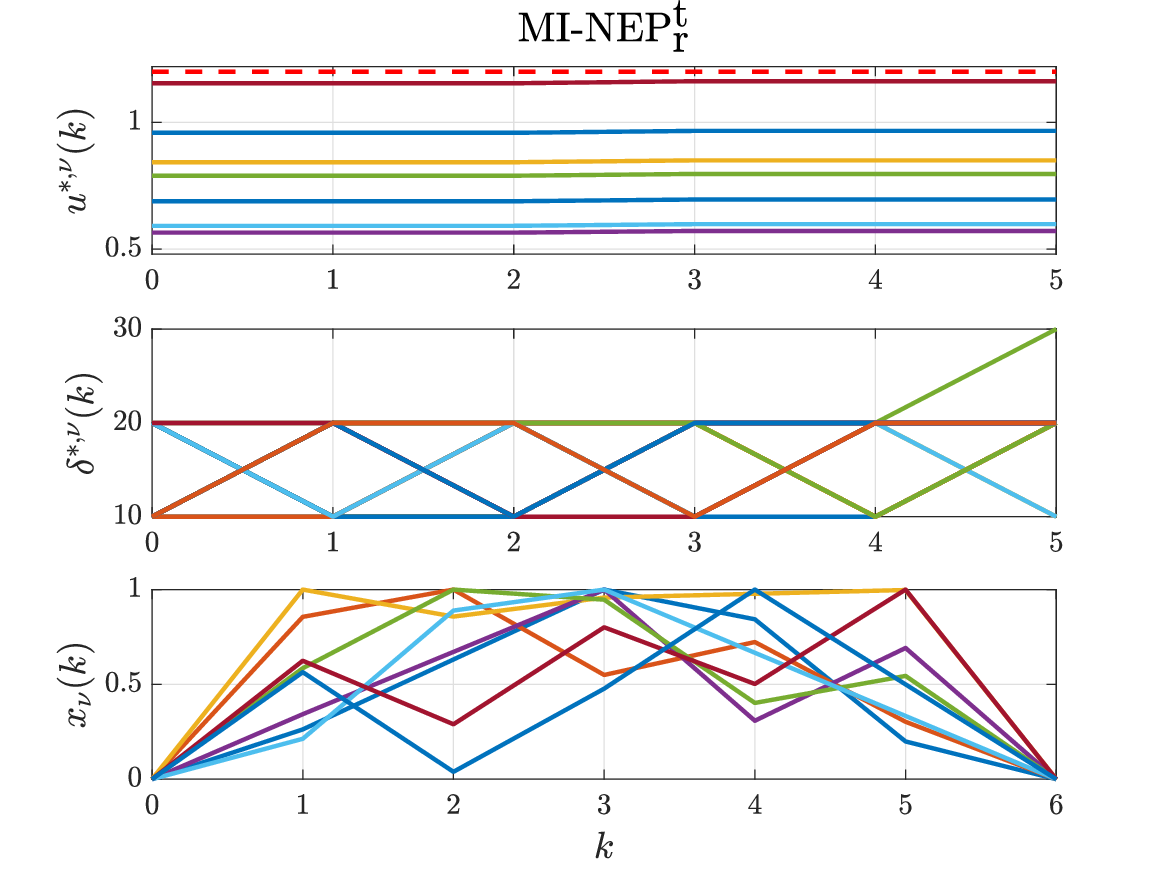}  
		\caption{}
		\label{fig:sub-fifth}
	\end{subfigure}
	\caption{Example of equilibria computed in all the considered cases. Specifically, (a) Relaxed \gls{NEP}; \gls{MI-NE} with granularity up to units over the entire feasible set (b) and over the reduced one (c); \gls{MI-NE} with granularity up to tens over the entire feasible set (d) and over the reduced one (e).}
	\label{fig:results}
\end{figure}

The plots in Fig.~\ref{fig:results}, instead, report five equilibria for a numerical instance in which Algorithm~\ref{alg: Jacobi} with $\mc J^k=\mc I$, for all $k$, has converged in all the considered cases \gls{MI-NEP}$^\triangle_\ast$, $\triangle \in \{\text{u}, \text{t}\}$, $\ast \in \{\text{f}, \text{r}\}$. Relaxing the integer restriction on $\delta^\nu(k)$ results in a continuous discharging of the battery, which generates a smoother \gls{SOC} $x_\nu(k)$ in Fig.~\ref{fig:results}(a) compared to those in Fig.~\ref{fig:results}(b)--(e). These latter are indeed strictly dependent on the integral restrictions on the scheduling variable $\delta^\nu$. In addition, from the evolution of the \gls{SOC} in Fig.~\ref{fig:results}(a) we can also appreciate the effect of the daily and night tariffs for energy consumption, which is not so evident with the \gls{MI-NE} in Fig.~\ref{fig:results}(b)--(e). While the equilibra computed for the \gls{MI-NEP}$^\text{u}_\ast$, $\ast \in \{\text{f}, \text{r}\}$, in Fig.~\ref{fig:results}(b)--(c) appear almost coinciding, with a scheduling variable $\delta^\nu$ assuming integer values between $15$ and $17$, the equilibria for the cases \gls{MI-NEP}$^\text{t}_\ast$, $\ast \in \{\text{f}, \text{r}\}$, in Fig.~\ref{fig:results}(d)--(e) instead are substatially different, mostly because $\delta^\nu$ varies between $0$ and $40$.

\section{Conclusion}
Focusing on traditional \gls{BR} algorithms, we have characterized the convergence properties of the sequence produced in computing solutions to a wide class of \glspl{MI-NEP}, i.e., problems that turn into monotone \glspl{NEP} once relaxed the integer restrictions. In particular, we have shown that the resulting sequence always approaches a bounded region containing the entire solution set of the \gls{MI-NEP}, whose size depends on the problem data. Moreover, we have confirmed that, once a Jacobi/Gauss-Seidel \gls{BR} method is applied to the relaxed \gls{NEP},  it converges to the unique solution, and we have also established data-dependent complexity results to characterize its convergence. Nonetheless, we have derived a sufficient condition for the existence of solutions to \glspl{MI-NEP}, as well as investigated the relation between the contraction property of the continuous \gls{NEP} and the degree of strong monotonicity possessed by such a relaxed problem. Numerical results on an instance of a smart building control application have illustrated the practical advantages brought by our results.

Future work will be devoted to analyze the numerical benefit entailed by the proposed results when combined with enumerative procedures, as well as to investigate the tight relation between (strong) monotonicity of a \gls{MI-NEP} and the resulting existence/uniqueness of a solution, thus fully characterizing the connection between Propositions~\ref{pr: strong monotonicity} and  \ref{pr: existence}.

\backmatter

\begin{appendices}

\section*{Appendix}

\begin{proposition}\label{pr: block norn}
 The function $\|\cdot\|^{\cal{B}(\wbold)}$, with $\wbold \geq 0$, is a seminorm, namely
 \begin{enumerate}[(i)]
  \item $\|\zbold + \ybold\|^{\cal{B}(\wbold)} \leq \|\zbold\|^{\cal{B}(\wbold)} + \|\ybold\|^{\cal{B}(\wbold)}$ for all $\zbold, \ybold \in \Re^n$,
  \item $\|a \zbold\|^{\cal{B}(\wbold)} = |a| \|\zbold\|^{\cal{B}(\wbold)}$ for all $\zbold \in \Re^n$ and $a \in \Re$.
 \end{enumerate}
 If $\wbold > 0$, then $\|\cdot\|^{\cal{B}(\wbold)}$ is a norm, i.e., in addition to (i) and (ii), it satisfies:
 \begin{enumerate}[(iii)]
  \item $\|\zbold\|^{\cal{B}(\wbold)} = 0$ implies $\zbold = 0$.
 \end{enumerate}
 \hfill$\square$
\end{proposition}
\textit{Proof.}
 (i) Let $\overline \nu \in \arg \max_{\nu \in \mc I}~w_\nu \|z^\nu + y^\nu\|$. We have: $\|\zbold + \ybold\|^{\cal{B}(\wbold)}  = w_{\overline \nu} \|z^{\overline \nu} + y^{\overline \nu}\| \leq w_{\overline \nu} (\|z^{\overline \nu}\| + \|y^{\overline \nu}\|) \leq \max_{\nu \in \mc I}~w_\nu (\|z^\nu\| + \|y^\nu\|) = \|\zbold\|^{\cal{B}(\wbold)} + \|\ybold\|^{\cal{B}(\wbold)}$.

\noindent
 (ii)
 $
 \|a \zbold\|^{\cal{B}(\wbold)} = \max_{\nu \in \mc I}~w_\nu \|a z^\nu\| = |a| \max_{\nu \in \mc I}~w_\nu \|z^\nu\| = |a| \|\zbold\|^{\cal{B}(\wbold)}.
 $
 
 \noindent
 (iii) Since $\wbold > 0$, $\|\zbold\|^{\cal{B}(\wbold)} = 0 \Rightarrow \|z^\nu\| = 0$, i.e., $z^\nu = 0$ for all $\nu \in \mc I$.
 \hfill$\blacksquare$
\end{appendices}


\bibliography{Surbib}

\end{document}